\documentclass{lms}
\usepackage[T1]{fontenc}
\usepackage[utf8]{inputenc}
\usepackage[english]{babel}
\usepackage{amsmath,amssymb,fancyhdr,enumerate}
\usepackage[all]{xy}
\title[Degree and representation type]{Degrees of irreducible
  morphisms and finite-representation type}
\author{Claudia Chaio, Patrick Le Meur, Sonia Trepode}
\date{\today}

\classno{16G60, 16G70}

\extraline{The first and third
 authors acknowledge partial support of CONICET, Argentina. The third
 author is a researcher of CONICET.
The third author acknowledges financial support from the Department of Mathematics of FCEyN, Universidad Nacional
 de Mar del Plata and also from the PICS 3410 of CNRS, France.}

\newtheorem{prop}{Proposition}[section]
\newtheorem{Thm}{Theorem}

\newtheorem{cor}[prop]{Corollary}
\newtheorem{lem}[prop]{Lemma}

\newnumbered{ex}[prop]{Example}
\newnumbered{rem}[prop]{Remark}
\newnumbered{defn}[prop]{Definition}

\def\ts#1{\normalfont{\textsf{#1}}}

\begin{document}

\maketitle
\abstract{We study the degree of irreducible morphisms in any
  Auslander-Reiten component of a finite dimensional algebra over an
  algebraically closed field. We give a characterization for an
  irreducible morphism to have finite left (or right) degree. This is
  used to prove our main theorem:  An algebra is of finite representation type if and only 
  if for every indecomposable projective the inclusion
  of the radical in the projective has finite right degree, which is
  equivalent to require that
  for every indecomposable injective the epimorphism from the
  injective to its quotient by its socle has finite
  left degree. We also apply the techniques that we develop: We study when
  the non-zero composite of a path of $n$ irreducible morphisms between
  indecomposable modules lies in the $n+1$-th power of the radical; and
  we study the same problem for sums of such paths when they are
  sectional, thus proving a generalisation of a pioneer result of
  Igusa and Todorov on the composite of a sectional path.}

\section*{Introduction}
Let $A$ be an artin algebra over an artin commutative ring $\ts k$. The
representation theory of $A$ deals with the study of the category
$\ts{mod}\,A$ of (right) $A$-modules of finite type. One of the most
powerful tools in this study is the Auslander-Reiten theory, based on
irreducible morphisms and almost split sequences (see
\cite{ARS95}). Although irreducible morphisms have permitted important
advances  in representation theory, some of their basic properties
still remain mysterious to us. An important example is the composition
of two irreducible morphisms: It obviously lies in $\ts{rad}^2$ (where
$\ts{rad}^l$ is the $l$-th power of the radical ideal $\ts{rad}$ of
$\ts{mod}\,A$) but it may lie in $\ts{rad}^3$, $\ts{rad}^{\infty}$ or
even be the zero morphism. Of course, the situation still makes sense
with the composite of arbitrary many irreducible morphisms. A first,
but partial, treatment of this situation was given by Igusa and
Todorov (\cite{IT84}) with the following result: ``\textit{If
  $X_0\xrightarrow{f_1}X_1\to\cdots\to X_{n-1}\xrightarrow{f_n}X_n$ is a
  sectional path of irreducible morphisms between indecomposable modules,
  then the composite $f_n\cdots f_1$ lies in $\ts{rad}^n(X_0,X_n)$
  and not in $\ts{rad}^{n+1}(X_0,X_n)$, in particular, it is
  non-zero.}'' In \cite{L92}, Liu introduced the left and right
degrees of an irreducible morphism $f\colon X\to Y$ as follows: The 
\emph{left degree} $d_l(f)$ of $f$ is the least integer $m\geqslant 1$
such that there exists $Z\in\ts{mod}\, A$ and
$g\in\ts{rad}^m(Z,X)\backslash\ts{rad}^{m+1}(Z,X)$ satisfying
$fg\in\ts{rad}^{m+2}(Z,Y)$. If no such an integer $m$ exists, then
$d_l(f)=\infty$. The \emph{right degree} is defined dually. This
notion was introduced to study the composition of irreducible
morphisms. In particular, Liu extended the above study of Igusa and
Todorov to presectional paths. Later it was used to determine the possible shapes of the
Auslander-Reiten components of $A$ (see \cite{L92,L96}). 
More recently, the composite of irreducible morphisms was studied in
\cite{CCT07}, \cite{CCT08}, \cite{CCT08a} and \cite{CT08}. The work
made in the first three of these papers is based on the notion of degree of
irreducible morphisms.
The definition
of the degree raises the following problem: Determine when $d_l(f)=\infty$ or
$d_r(f)=\infty$. Consider an irreducible morphism $f\colon X\to Y$
with $X$ indecomposable. Then, the following conditions have been
related in the recent literature:
\begin{enumerate}
\item $d_l(f)=n<\infty$,
\item $\ts{Ker}(f)$ lies in the Auslander-Reiten component containing
  $X$.
\end{enumerate}
Indeed, these two conditions were proved to be equivalent if the
Auslander-Reiten component containing $X$ is convex,
generalized standard and with length (\cite{CPT04}, actually this
equivalence still holds true if one removes the convex hypothesis) and when the
Auslander-Reiten quiver is standard
 (\cite{C08}).
In this text, we shall see that such results are key-steps to
show that the degree of irreducible morphisms is a useful notion to
determine the representation type of $A$. Indeed, we recall the
following well-known conjecture appeared first in \cite{L96} and
related to the Brauer-Thrall conjectures: "\textit{If the Auslander-Reiten
  quiver of $A$ is connected, then $A$ is of finite representation
  type.}"  This conjecture is related to the degree of irreducible
morphisms as follows: In the above situation of assertions (1) and
(2), the existence of $f$ such that $d_l(f)=\infty$ is related to the
existence of at least two Auslander-Reiten components. 
Actually, it was proved in \cite[Thm. 3.11]{CPT04} that if $A$ is of
finite representation type, then every irreducible morphism between
indecomposables either has finite right degree or has finite left degree.
Conversely, one can wonder if the converse holds true.
In this text, we prove the following main theorem where we assume that
$\ts k$ is an algebraically closed field.
\begin{Thm}
\label{thm1}
  Let $A$ be a connected finite dimensional $\ts k$-algebra over an
  algebraically closed field. The following conditions are equivalent:
  \begin{enumerate}[(a)]
\item $A$ is of finite representation type.
  \item For every indecomposable projective $A$-module $P$,
    the inclusion $\ts{rad}(P)\hookrightarrow P$ has finite right
    degree.
\item For every indecomposable injective $A$-module $I$, the quotient $I\to
  I/\ts{soc}(I)$ has finite left degree. 
\item For every irreducible epimorphism $f\colon X\to Y$ with $X$ or
  $Y$ indecomposable, the left degree of $f$ is finite.
\item For every irreducible monomorphism $f\colon X\to Y$ with $X$ or
  $Y$ indecomposable, the right degree of $f$ is finite.
  \end{enumerate}
\end{Thm}
Hence, going back to the above conjecture, if one knows that the
Auslander-Reiten quiver of $A$ is connected, by (b) and (c) it
suffices to study the
degree of finitely many irreducible morphisms in order to prove that
$A$ is of finite representation type.
Our proof of the above theorem only uses considerations on degrees and
their interaction 
with coverings of translation quivers. In particular it uses no
advanced characterization of finite representation type (such as the
Brauer-Thrall conjectures or multiplicative bases, for example). The
theorem shows that the degrees of irreducible
morphisms are somehow related to the representation type of $A$. Note
also that our characterization is expressed in terms of the knowledge
of the degree of finitely many irreducible morphisms. In order to
prove the theorem we investigate the degree of irreducible morphisms
and more particularly assertions (1) and (2) above. Assuming that $\ts k$
is an algebraically closed field and given $f\colon X\to Y$ an
irreducible epimorphism with $X$ indecomposable,  we prove that the
assertion (1) is equivalent to (3) below and implies (2), with no
assumption on the Auslander-Reiten component $\Gamma$ containing $X$:
\begin{enumerate}
\setcounter{enumi}{2}
\item There exists $Z\in\Gamma$ and
  $h\in\ts{rad}^n(Z,X)\backslash\ts{rad}^{n+1}(Z,X)$ such that $fh=0$.
\end{enumerate}
Therefore, the existence of an irreducible monomorphism (or
epimorphism) with infinite left (or right)
degree indicates  that
there are more than one  component in the Auslander-Reiten quiver (at
least when $\Gamma$ is generalized standard).
We also prove that (2) implies (1) (and therefore implies (3)) under the
additional assumption that $\Gamma$ is generalized standard. The
equivalence between (1) and (3) and the fact that it works for any
Auslander-Reiten component are the chore facts in the proof of the
theorem. For this purpose we use the covering techniques introduced
in \cite{R80a}. Indeed, these techniques allow one to reduce the study of
the degree of irreducible morphisms in a component to the study of the degree of
irreducible morphisms in a suitable covering called the generic
covering. Among other things, the
generic covering  is a translation quiver with length. As was proved in
\cite{CPT04} such a condition is particularly useful in the study of
the degree of an irreducible morphism. 

The text is therefore organized as follows. In the first section we
recall some needed definitions. In the second section we extend to any
Auslander-Reiten component the pioneer result \cite[2.2, 2.3]{R80a} on
covering techniques
which, in its original form, only works for 
the Auslander-Reiten quiver of representation-finite algebras. The
results of this section are used in the third one to prove the
various implications between assertions (1), (2) and (3) in
\ref{prop:degree} and \ref{prop:kernel}. As explained above, these
results have been studied previously and they were proved under
additional assumptions. In particular, the
corresponding corollaries proved at that time can be generalized
accordingly. In the fourth
section we prove our main Theorem~\ref{thm1} using the previous
results. The proof of our main results are based on the covering
techniques developed in the second section. In the last section, we
use these to study when the non-zero
composite
of $n$ irreducible morphisms lies in the $n+1$-th power of the
radical and we extend the cited-above result (\cite{IT84}) of Igusa and Todorov on
the composite of a sectional paths
to sums of composites of sectional paths.

\section{Preliminaries}

\subsection{Notations on modules}

Let $A$ be a finite dimensional $\ts k$-algebra. We denote by $\ts{ind}\, A$ a
full subcategory of $\ts{mod}\, A$ which contains exactly one representative
of each isomorphism class of indecomposable modules. Also, we write
$\ts{rad}$ for the \emph{radical} of $\ts{mod}\, A$. Hence, given indecomposable
modules $X,Y$, the space $\ts{rad}(X,Y)$ is the subspace of $\ts{Hom}_A(X,Y)$
consisting of non-isomorphisms $X\to Y$. For $l\geqslant 1$, we write
$\ts{rad}^l$ for the $l$-th power of the ideal $\ts{rad}$, recursively defined
by $\ts{rad}^1=\ts{rad}$ and $\ts{rad}^{l+1}=\ts{rad}.\ts{rad}^l$
($=\ts{rad}^l.\ts{rad}$). 
For short we shall say that some morphisms $u_1,\ldots,u_r\colon X\to
Y$ are \emph{linearly independent modulo $\ts{rad}^n(X,Y)$} if their
respective classes modulo $\ts{rad}^n(X,Y)$ are linearly independent
in $\ts{Hom}_A(X,Y)/\ts{rad}^n(X,Y)$.
We
recall that the \emph{Auslander-Reiten quiver of $A$} is the
translation quiver $\Gamma(\ts{mod}\, A)$ with vertices the modules in $\ts{ind}\, A$,
such that the number of arrows $X\to Y$ equals the dimension of the
quotient space $\ts{rad}(X,Y)/\ts{rad}^2(X,Y)$ for every vertices $X,Y\in\Gamma$
and whose translation is induced by the Auslander-Reiten translation
$\tau_A=D\ts{Tr}$. Hence, the translation quivers we shall deal with
are not valued quivers and may have multiple (parallel) arrows. If
$\Gamma$ is a component of $\Gamma(\ts{mod}\,
A)$ (or an \emph{Auslander-Reiten component}, for short), we write $\ts{ind}\,\Gamma$
for the full subcategory of $\ts{ind}\, A$ with objects the modules in
$\Gamma$.
Recall that a \emph{hook} is a path $X\to Y\to Z$ of
irreducible morphisms between indecomposable modules such that $Z$ is
non-projective and $X=\tau_AZ$. Also, a path
$X_0\to X_1\to\cdots\to X_{l-1}\to X_l$ of irreducible morphisms is \emph{sectional} if neither
of its subpaths of length $2$ is a hook.

We refer the reader to \cite{L92} for properties on the degree
of irreducible morphisms.

\subsection{Radical in mesh-categories}

Let $\Gamma$ be a \emph{translation quiver}, that is, $\Gamma$ is a quiver
with no loops (but possibly with parallel arrows); endowed with two
distinguished subsets of vertices, the elements of
which are called projectives and injectives, respectively; and endowed
with a bijection $\tau\colon x\mapsto \tau x$ (the \emph{translation})  from the set of
non-projectives to the set of non-injectives; such that for every
vertices $x,y$ with $x$ non-projective, there is a bijection
$\alpha\mapsto \sigma\alpha$ from the set of arrows $y\to x$ to the
set of arrows $\tau x\to y$. 
All translation quivers are assumed to be \emph{locally finite}, that
is, every vertex is the source or the target of at most finitely many
arrows (Auslander-Reiten components are always locally finite quivers).
The subquiver of $\Gamma$ formed by the
arrows starting at $\tau x$ and the arrows arriving at $x$ is called
the \emph{mesh} ending at $x$. We write $\ts k(\Gamma)$ for the
\emph{mesh-category} of $\Gamma$, that is, the factor category of the
path category $\ts k\Gamma$ by the ideal generated by the morphisms
$\sum\limits_{\alpha\colon \cdot\to x}\alpha\ \sigma\alpha$ where $\alpha$
runs through the arrows arriving at $x$, for a given non-projective
vertex $x$. If $u$ is a path in $\Gamma$, we write
$\overline u$ for the corresponding morphism in $\ts k(\Gamma)$. We denote
by  $\mathfrak{R}\ts k(\Gamma)$ the ideal of
$\ts k(\Gamma)$ generated by $\{\overline
\alpha\ |\ \text{$\alpha$ an arrow in $\Gamma$}\}$. 
Note that in general $\mathfrak{R}\ts k(\Gamma)$ is not a radical of the
category $\ts k(\Gamma)$.
The $l$-th power
$\mathfrak{R}^l\ts k(\Gamma)$ is defined recursively by $\mathfrak{R}^1\ts k(\Gamma)=\mathfrak{R}\ts k(\Gamma)$ and
$\mathfrak{R}^{l+1}\ts k(\Gamma)=\mathfrak{R}\ts k(\Gamma).\mathfrak{R}^l\ts k(\Gamma)$
($=\mathfrak{R}^l\ts k(\Gamma).\mathfrak{R}\ts k(\Gamma)$). If $\Gamma$
is
\emph{with length}, that is, for every $x,y\in\Gamma$ all the paths from
$x$ to $y$ have equal length, then the radical satisfies the following
result proved in \cite[Prop. 2.1]{CT08}. This result is central in our
work. Later we shall use it without further reference.
\begin{prop}
\label{prop:with_length}
  Let $\Gamma$ be a translation quiver with length and $x,y\in\Gamma$. If
  there is a path of length $l$ from $x$ to $y$ in $\Gamma$, then:
  \begin{enumerate}[(a)]
  \item $\ts k(\Gamma)(x,y)=\mathfrak{R}\ts k(\Gamma)(x,y)=\mathfrak{R}^2\ts k(\Gamma)(x,y)=\ldots=\mathfrak{R}^l\ts k(\Gamma)(x,y)$.
  \item $\mathfrak{R}^i\ts k(\Gamma)(x,y)=0$ if $i>l$.
  \end{enumerate}
\end{prop}
In view of the preceding proposition, we call a \emph{length function
  on $\Gamma$} a function $l$ which assigns an integer
$l(x)\in\mathbb{Z}$ to every vertex $x\in\Gamma$, in such a way that
$l(y)=l(x)+1$ whenever there is an arrow $x\to y$ in
$\Gamma$ (see \cite[1.6]{BG82}). Clearly, if $\Gamma$ has a length function, then $\Gamma$ is
with length.
Finally, we define hooks and sectional paths in
translation quivers as we did for hooks and sectional paths of
irreducible morphisms in module categories.

\subsection{Coverings of translation quivers}

Let $\Gamma$ be a connected translation quiver. A \emph{covering of
  translation quivers} (\cite[1.3]{BG82}) is a morphism $p\colon \Gamma'\to \Gamma$ of quivers
such that:
\begin{enumerate}[(a)]
\item $\Gamma'$ is a translation quiver.
\item A vertex $x\in\Gamma'$ is projective (or injective,
  respectively) if and only if so is $px$.
\item $p$ commutes with the translations in $\Gamma$ and $\Gamma'$
  (when these are defined).
\item For every vertex $x\in\Gamma'$, the map $\alpha\mapsto
  p(\alpha)$ induces a bijection from the set of arrows in $\Gamma'$
  starting from $x$ (or ending at $x$) the set of arrows in $\Gamma$
  starting from $p(x)$ (or ending at $p(x)$, respectively).
\end{enumerate}
We shall use a particular covering
$\pi\colon\widetilde{\Gamma}\to\Gamma$ which we call the \emph{generic
covering}. Following \cite[1.2]{BG82}, we define the equivalence
relation $\sim$ on the set of unoriented paths in $\Gamma$ as
generated by the following properties
\begin{enumerate}[(i)]
\item If $\alpha\colon x\to y$ is an arrow in $\Gamma$, then
  $\alpha\alpha^{-1}\sim e_y$ and $\alpha^{-1}\alpha\sim e_x$ (where
  $e_x$ denotes the stationary path at $x$, of length $0$).
\item If $x$ is a non-projective vertex and the mesh in $\Gamma$
  ending at $x$ has the form
  \begin{equation}
    \xymatrix@R=1ex{
&x_1 \ar[rd]^{\beta_1} \ar@{.}[dd]\\
\tau x \ar[ru]^{\alpha_1} \ar[rd]_{\alpha_r} && x
\\
&x_r \ar[ru]_{\beta_r}
}\notag
  \end{equation}
then $\beta_i\alpha_i\sim\beta_j\alpha_j$ for every
$i,j\in\{1,\ldots,r\}$.
\item If $\alpha,\beta$ are arrows in $\Gamma$ with the same source
  and the same target, then $\alpha\sim\beta$.
\item If $\gamma_1,\gamma,\gamma',\gamma_2$ are unoriented paths such
  that $\gamma\sim\gamma'$ and the compositions
  $\gamma_1\gamma\gamma_2,\gamma_1\gamma'\gamma_2$ are defined, then
  $\gamma_1\gamma\gamma_2\sim\gamma_1\gamma'\gamma_2$.
\end{enumerate}
Note that the usual homotopy relation of a translation quiver
(see \cite[1.2]{BG82}) is defined using conditions ($i$),
($ii$) and ($iv$) above. Also recall that the universal cover of
$\Gamma$ was defined in \cite[1.3]{BG82} using that homotopy
relation. By applying that construction to our equivalence relation
$\sim$ instead of to the homotopy relation, we get the covering
$\pi\colon \widetilde{\Gamma}\to \Gamma$ which we call the
\emph{generic covering} of $\Gamma$. 
Note that if $\Gamma$ has no multiple arrows (for example, if $\Gamma=\Gamma(\ts{mod}\,A)$ is the Auslander-Reiten
quiver of a representation-finite algebra), then the generic covering
coincides with the universal covering.
The following properties of $\pi$ are crucial to our work:
\begin{prop}
\label{prop:generic_covering} Let $\Gamma$ be a translation quiver and
$\pi\colon\widetilde{\Gamma}\to \Gamma$ be its generic covering.
\begin{enumerate}[(a)]
\item There is a length function on $\widetilde{\Gamma}$. In
  particular, $\widetilde{\Gamma}$ is with length.
\item If $\alpha\colon x\to y$, $\beta\colon x\to z$ (or $\alpha\colon
  y\to x$, $\beta\colon z\to x$) are arrows in $\widetilde{\Gamma}$
  such that $\pi y=\pi z$, then $y=z$.
\item For every vertices $x,y\in\widetilde{\Gamma}$, the covering
  $\pi$ induces a bijection from the set of arrows in
  $\widetilde{\Gamma}$ from $x$ to $y$ to the set of arrows in
  $\Gamma$ from $\pi x$ to $\pi y$.
\item Let $x,y\in\widetilde{\Gamma}$ be vertices. If $u\colon x=x_0\to
  \cdots\to x_l=y$ and $v\colon x=x_0'\to\cdots\to x_l'=y$ are
  two paths in $\widetilde{\Gamma}$ from $x$ to $y$ and if $u$ is
  sectional, then $x_1=x_1',\ldots,x_{l-1}=x_{l-1}'$. In particular,
  all the paths from $x$ to $y$ are sectional.
\item Let $x,y\in\widetilde{\Gamma}$ be vertices, $u_1,\ldots,u_r$
  be pairwise distinct sectional paths of length $n\geqslant 0$ in
  $\widetilde{\Gamma}$ from $x$ to $y$ and
  $\lambda_1,\ldots,\lambda_r\in \ts k$ be scalars. Then the following
  equivalence holds in $\ts k(\widetilde{\Gamma})$:
  \begin{equation}
    \lambda_1\overline{u_1}+\cdots+\lambda_r\overline{u_r}\in
    \mathfrak{R}^{n+1}\ts k(\widetilde{\Gamma}) \Leftrightarrow \lambda_1=\cdots=\lambda_r=0.\notag
  \end{equation}
\end{enumerate}
\end{prop}
\begin{proof}
  Let $\Gamma'$ be the translation quiver with no multiple arrows, with the same vertices and
  the same translation as those of $\Gamma$, and such that there is an
  arrow (and exactly one) $x\to y$ in $\Gamma'$ if and only if there
  is (at least) one arrow $x\to y$ in $\Gamma$. Define
  $\widetilde{\Gamma'}$ starting from $\widetilde{\Gamma}$ in a
  similar way. Then $\widetilde{\Gamma'}$ is the universal cover of
  $\Gamma'$ and it is simply connected (in the sense of \cite[1.3,1.6]{BG82}).\\

(a) Applying \cite[1.6]{BG82} to $\widetilde{\Gamma'}$ yields a length
function on $\widetilde{\Gamma'}$ and, therefore, on
$\widetilde{\Gamma}$. The existence of a length function implies that
$\widetilde{\Gamma}$ is with length.

(b) follows from the construction of the generic covering.

(c) follows from (b) and from the fact that $\pi\colon
\widetilde{\Gamma}\to\Gamma$ is a covering of quivers.

(d) The paths $u$ and $v$ define two (unique) paths $x_0\to
  \cdots\to x_l$ and $x_1'\to\cdots\to x_l'$ from $x$ to $y$ in
  $\widetilde{\Gamma'}$, the first of which is sectional. The
  conclusion then follows from \cite[Lem. 1.2]{BL03}.

(e) Assume that
$\lambda_1\overline{u_1}+\cdots+\lambda_r\overline{u_r}\in\mathfrak{R}^{n+1}\ts k(\widetilde{\Gamma})$. It
follows from (a) and from \ref{prop:with_length} that
$\lambda_1\overline{u_1}+\cdots+\lambda_r\overline{u_r}=0$, that is,
$\lambda_1 u_1+\cdots+\lambda_r u_r$ lies in the mesh-ideal. By
definition of the mesh-ideal, this implies that $u_i$ contains a hook
whenever $\lambda_i\neq 0$. Using (d), we deduce that
$\lambda_1=\cdots=\lambda_r=0$. The converse is obvious.
\end{proof}

The second property (b) in \ref{prop:generic_covering} is not satisfied by the universal cover
when $\Gamma$ has multiple arrows. This is the reason for using the
generic covering instead.

\section{Well-behaved functors}
\label{sec:s1}
Let $A$ be a finite dimensional $\ts k$-algebra where $\ts k$ is an
algebraically closed field, $\Gamma$  a
component of $\Gamma(\ts{mod}\, A)$ and
$\pi\colon\widetilde{\Gamma}\to\Gamma$ the generic covering. Following
\cite[3.1 Ex. (b)]{BG82} (see also \cite{R80a}), a
$\ts k$-linear functor $F\colon \ts k(\widetilde{\Gamma})\to
\ts{ind}\,\Gamma$ is called \emph{well-behaved} if it satisfies the
following conditions for every vertex $x\in\widetilde{\Gamma}$: (a)
$Fx=\pi x$; (b) If
$\alpha_1\colon x\to x_1,\ldots,\alpha_r\colon x\to x_r$ are the
arrows in $\widetilde{\Gamma}$ starting from $x$, then $
[F(\overline\alpha_1),\ldots , F(\overline\alpha_r)]^t
\colon Fx\to Fx_1\oplus\cdots \oplus Fx_r$ is
minimal left almost split; (c) 
If
$\alpha_1\colon  x_1\to x,\ldots,\alpha_r\colon x_r\to x$ are the
arrows in $\widetilde{\Gamma}$ ending at $x$, then $
[F(\overline\alpha_1),\ldots , F(\overline\alpha_r)]
\colon  Fx_1\oplus\cdots \oplus Fx_r\to Fx$ is
minimal right almost split. Note that these conditions imply that $F$
maps meshes in
$\widetilde{\Gamma}$ to almost split sequences in $\ts{mod}\,A$.

For convenience, we extend this notion to functors $p\colon 
\ts k\mathcal X\to\ts{ind}\,\Gamma$ where $\mathcal X$ is a subquiver of
$\widetilde{\Gamma}$. The functor $p$ is called well-behaved if and
only if ($1$) $px=\pi x$ for every vertex $x\in\mathcal X$; ($2$)
given a vertex $x\in\mathcal X$, if
$x\xrightarrow{\alpha_1}x_1,\ldots, x\xrightarrow{\alpha_r}x_r$ are
the arrows in $\mathcal X$ starting in $x$, then the morphism
$[p(\alpha_1),\ldots,p(\alpha_r)]^t\colon \pi x\to
\bigoplus\limits_{i=1}^r\pi x_i$ is irreducible; ($3$)
given a vertex $x\in\mathcal X$, if
$x_1\xrightarrow{\beta_1}x,\ldots, x_r\xrightarrow{\beta_r}x$ are
the arrows in $\mathcal X$ ending in $x$, then the morphism
$[p(\beta_1),\ldots,p(\beta_r)]\colon 
\bigoplus\limits_{i=1}^r\pi x_i\to \pi x$ is irreducible; ($4$) if the
vertex $x$ is non-projective and if $\mathcal X$ contains the mesh in
$\widetilde{\Gamma}$ ending in $x$
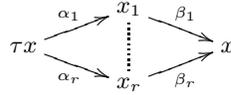
  \begin{equation}
    \xymatrix@R=1ex{
&x_1\ar[rd]^{\beta_1} \ar@{.}[dd]&\\
\tau x\ar[ru]^{\alpha_1}\ar[rd]_{\alpha_r} && x\\
&x_r \ar[ru]_{\beta_r}
}\notag
  \end{equation}
then the sequence $0\to \tau_A\pi
x\xrightarrow{[p(\alpha_1),\ldots,p(\alpha_r)]^t}\bigoplus\limits_{i=1}^r \pi x_i
\xrightarrow{[p(\beta_1),\ldots,p(\beta_r)]} \pi
x\to 0$ is exact and almost split. Recall that if
$X\xrightarrow{a_1}X_1,\ldots, X\xrightarrow{a_r}X_r$ are all the
arrows in $\Gamma$ starting in some module $X$ and if the morphism $[u_1,\ldots,
u_r]^t\colon X\to \bigoplus\limits_{i=1}^r X_i$ is irreducible, then
it is minimal left almost split (and dually). Therefore, a
well-behaved functor $p\colon \ts k\widetilde{\Gamma}\to \ts{ind}\,\Gamma$
induces a well-behaved functor $F\colon \ts
k(\widetilde{\Gamma})\to\ts{ind}\,\Gamma$ by factoring out by the
mesh-ideal.

Recall that in the case where $A$ is of finite representation type, it
was proved in \cite[\S 3]{BG82} that there always exists a
well-behaved functor $\ts
k(\widetilde{\Gamma})\to\ts{ind}\,\Gamma$. This result was based on a
similar one in \cite[\S 1]{R80a}, where a well-behaved functor $\ts
k(\widetilde{\Gamma})\to \underline{\ts{mod}}\,A$ was constructed when
$A$ is self-injective and of finite representation type and $\Gamma$
is a stable component of $\Gamma(\ts{mod}\,A)$.

In this text we use the following more general existence result on
well-behaved functors.  Given a length
 function $l$ on $\widetilde{\Gamma}$ and an integer $n\in\mathbb{Z}$, we denote by
 $\widetilde{\Gamma}_{\leqslant n}$ (or
 $\widetilde{\Gamma}_{\geqslant n}$) the full subquiver of
 $\widetilde{\Gamma}$ with vertices those $x\in\widetilde{\Gamma}$
 such that $l(x)\leqslant n$ (or $l(x)\geqslant n$,
 respectively). These are convex subquivers.
 \begin{prop}
   \label{prop:well_behaved_general}
Let $\pi\colon \widetilde{\Gamma}\to \Gamma$ be the generic covering
and let $q\colon  \ts k\mathcal Y\to\ts{ind}\,\Gamma$ be a well-behaved
functor with $\mathcal Y$ a full convex subquiver of
$\widetilde{\Gamma}$. Let $l$ be a length function on
$\widetilde{\Gamma}$ and assume that at least one of the following
conditions is satisfied:
\begin{enumerate}[(a)]
\item There exist integers $m,n\in\mathbb{Z}$ such that $n\leqslant m$
  and  $\mathcal
  Y\subseteq \widetilde{\Gamma}_{\geqslant n}\cap
  \widetilde{\Gamma}_{\leqslant m}$.
\item There exists an integer $n\in\mathbb{Z}$ such that $\mathcal
  Y\subseteq \widetilde{\Gamma}_{\geqslant n}$ and $\mathcal Y$ is
  stable under predecessors in $\widetilde{\Gamma}_{\geqslant n}$
  (that is, every path in $\widetilde{\Gamma}_{\geqslant n}$ with
  endpoint lying in $\mathcal Y$ lies entirely in $\mathcal Y$).
\item There exists an integer $n\in \mathbb{Z}$ such that
  $\widetilde{\Gamma}_{\geqslant n}\subseteq\mathcal Y$ and $\mathcal
  Y$ is stable under successors in $\widetilde{\Gamma}$.
\end{enumerate}
Then there exists a well-behaved functor $F\colon \ts
k(\widetilde{\Gamma})\to \ts{ind}\,\Gamma$ such that
$F(\overline{\alpha})=q(\alpha)$ for every arrow $\alpha\in\mathcal Y$.
 \end{prop}
\begin{proof} It is sufficient to prove that there exists a
well-behaved functor $p\colon \ts k \widetilde{\Gamma}\to
\ts{ind}\,\Gamma$. For that purpose, we shall prove that:
\begin{itemize}
\item[-] If $\mathcal Y$ satisfies (a), then $q$ extends to a
  well-behaved functor $p\colon  \ts k\mathcal X\to \ts{ind}\,\Gamma$
  with $\mathcal X$ a full convex subquiver of $\widetilde{\Gamma}$
  satisfying (b).
\item[-] If $\mathcal Y$ satisfies (b), then $q$ extends to a
  well-behaved functor $p\colon \ts k\mathcal X\to\ts{ind}\,\Gamma$
  with $\mathcal X$ a full convex subquiver of $\widetilde{\Gamma}$
  satisfying (c).
\item[-] If $\mathcal Y$ satisfies (c), then $q$ extends to a
  well-behaved functor $p\colon \ts k\widetilde{\Gamma}\to \ts{ind}\,\Gamma$.
\end{itemize}
We shall consider pairs, $(\mathcal X,p)$ where $\mathcal X$ is a full
convex subquiver of $\widetilde{\Gamma}$ containing $\mathcal Y$ and $p\colon \ts k\mathcal
X\to \ts{ind}\,\Gamma$ is a well-behaved functor extending $q$. For
any two such pairs $(\mathcal X,p)$ and $(\mathcal X',p')$, we shall
write $(\mathcal X,p)\leqslant (\mathcal X',p')$ if and only if
$\mathcal X\subseteq\mathcal X'$ and $p'$ extends $p$. This clearly
defines a partial order on the set of such pairs.

Assume that $\mathcal Y$ satisfies (a). Consider the set $\Sigma$ of
those pairs $(\mathcal X,p)$ where $\mathcal X$ is a full convex
subquiver of $\widetilde{\Gamma}$ containing $\mathcal Y$, contained
in $\widetilde{\Gamma}_{\geqslant n}\cap \widetilde{\Gamma}_{\leqslant
  m}$, and $p\colon \ts k\mathcal X\to \ts{ind}\,\Gamma$ is a
well-behaved functor extending $q$. Then $\Sigma$ is non-empty since it
contains ($\mathcal X,p)$. Moreover, $(\Sigma, \leqslant)$ is totally
inductive. Therefore, it has a maximal element, say $(\mathcal
X,p)$. We claim that $\mathcal X$ is stable under predecessors in
$\widetilde{\Gamma}_{\geqslant n}$. By absurd, assume that this is not
the case. Then there exists an arrow $x\to y$ in
$\widetilde{\Gamma}_{\geqslant n}$ with $x\not\in\mathcal X$ and
$y\in\mathcal X$. We choose such an $x$ with $l(x)$ maximal. This is
possible because $\mathcal X\subseteq \widetilde{\Gamma}_{\leqslant
  m}$. Note that there is no arrow $z\to x$ in $\widetilde{\Gamma}$
with $z\in\mathcal X$, because, otherwise, the path $z\to x\to y$ would
contradict the convexity of $\mathcal X$. Therefore, the full
subquiver $\mathcal X'$ of $\widetilde{\Gamma}$ generated by $\mathcal
X$ and $x$ has, as arrows, those in $\mathcal X$ together with the
arrows in $\widetilde{\Gamma}$ starting in $x$ and ending at some
vertex in $\mathcal X$, say
\begin{equation}
  \xymatrix@R=1ex{
&x_1=y\ar@{.}[dd]\\
x \ar[ru]^{\alpha_1} \ar[rd]_{\alpha_r}\\
&x_r.
}\notag
\end{equation}
In particular, the convexity of $\mathcal X$ and the maximality of
$l(x)$ imply that $\mathcal
X'$ is convex. Assume that $x$ is injective, or else $x$ is
non-injective and $\tau^{-1}x\not\in\mathcal X$. The arrows
$\pi(\alpha_1),\ldots,\pi(\alpha_r)$ in $\Gamma$ are pairwise distinct
and start in $\pi x$. Therefore, there exists an irreducible morphism
$[u_1,\ldots, u_r]^t\colon \pi x\to \bigoplus\limits_{i=1}^r\pi
x_i$. We thus extend $p\colon \ts k\mathcal X\to \ts{ind}\,\Gamma$ to
a functor $p'\colon \ts k\mathcal X'\to \ts{ind}\,\Gamma$ by setting
$p'(\alpha_i)=u_i$, for every $i$. Note that a mesh in
$\widetilde{\Gamma}$ is contained in $\mathcal X$ if and only if it is
contained in $\mathcal X'$, by assumption on $x$ and because there is
no arrow in $\widetilde{\Gamma}$ ending at $x$ and starting in some
vertex in $\mathcal X$. Assume now that $x$ is
non-injective and $\tau^{-1}x\in\mathcal X$. By maximality of $l(x)$,
every arrow in $\widetilde{\Gamma}$ ending at $\tau^{-1}x$ lies in
$\mathcal X$. Therefore, $\alpha_1,\ldots,\alpha_r$ are all the arrows
in $\widetilde{\Gamma}$ starting in $x$ and the mesh in
$\widetilde{\Gamma}$ starting in
$x$ is as follows
\begin{equation}
  \xymatrix@R=1ex{
& x_1\ar[rd]^{\beta_1}\ar@{.}[dd]\\
x \ar[ru]^{\alpha_1} \ar[rd]_{\alpha_r} & & \tau^{-1}x\, .\\
& x_r \ar[ru]_{\beta_r}
}\notag
\end{equation}
Since $p$ is well-behaved, the morphism
$[p(\beta_1),\ldots,p(\beta_r)]\colon \bigoplus\limits_{i=1}^r \pi
x_i\to \tau_A^{-1}\pi x$ is irreducible and, therefore, minimal right
almost split, because $\pi(\beta_1),\ldots,\pi(\beta_r)$ are all the
arrows in $\Gamma$ ending in $\tau_A^{-1}\pi x$. Hence, there is an
almost split sequence in $\ts{mod}\,A$
\begin{equation}
  0\to \pi x\xrightarrow{[u_1,\ldots,u_r]^t} \bigoplus\limits_{i=1}^r
  \pi x_i \xrightarrow{[p(\beta_1),\ldots,p(\beta_r)]} \tau_A^{-1}\pi
    x\to 0\notag
\end{equation}
and we extend $p\colon \ts k\mathcal X\to \ts{ind}\,\Gamma$ to a
functor $p'\colon \ts k\mathcal X'\to\ts{ind}\,\Gamma$ by setting
$p(\alpha_i)=u_i$, for every $i$. In any case, $p'$ is
well-behaved. Indeed, by construction of $p'$ and because $p$ is
well-behaved, we have:  $p'$ transforms every mesh in $\widetilde{\Gamma}$
contained in $\mathcal X'$ into an almost split sequence in
$\ts{mod}\,A$; moreover, given vertices $z,t\in\mathcal X'$, if
$\gamma_1,\ldots,\gamma_s\colon z\to t$ are all the arrows in
$\mathcal X'$ from $z$ to $t$, then the morphism
$[p'(\gamma_1),\ldots,p'(\gamma_s)]^t\colon \pi z\to
\bigoplus\limits_{i=1}^s\pi t$ is irreducible (or, equivalently, the
morphism $[p'(\gamma_1),\ldots,p'(\gamma_s)]\colon \bigoplus\limits_{i=1}^s\pi z\to
\pi t$ is irreducible, because $\ts k$ is an algebraically closed field);
using \ref{prop:generic_covering}, (b), we deduce 
  that $p'$ is well-behaved. Thus $(\mathcal
X',p')\in\Sigma$ and $(\mathcal X,p)<(\mathcal X',p')$, a contradiction
to the maximality of $(\mathcal X,p)$. This proves that $q\colon \ts k
\mathcal Y\to \ts{ind}\,\Gamma$ extends to a well-behaved functor
$p\colon \ts k \mathcal X\to \ts{ind}\,\Gamma$ with $\mathcal X$ a
full convex subquiver of $\widetilde{\Gamma}$ containing $\mathcal Y$,
contained in $\widetilde{\Gamma}_{\geqslant n}\cap
\widetilde{\Gamma}_{\leqslant m}$ and stable under
predecessors in $\widetilde{\Gamma}_{\geqslant n}$. Therefore,
$\mathcal X$ satisfies (b).\\

Now assume that $\mathcal Y$ satisfies (b). Let $\Sigma'$ be the set
of those pairs $(\mathcal X,p)$ where $\mathcal X$ is a full convex
subquiver of $\widetilde{\Gamma}$ containing $\mathcal Y$, contained
in $\widetilde{\Gamma}_{\geqslant n}$  and stable under predecessors
in $\widetilde{\Gamma}_{\geqslant n}$, and $p\colon \ts k\mathcal X\to
\ts{ind}\,\Gamma$ is a well-behaved functor extending $q$. Then
$\Sigma'$ is non-empty for it contains $(\mathcal Y,q)$, and
$(\Sigma',\leqslant)$ is totally inductive. Let $(\mathcal X,p)$ be a
maximal element in $\Sigma'$. We claim that $\mathcal
X=\widetilde{\Gamma}_{\geqslant n}$. By absurd, assume that this
is not the case. Let $x\in\widetilde{\Gamma}_{\geqslant n}$ be a
vertex not in $\mathcal X$. We may assume that  $l(x)$ is
minimal for this property. Then $x$ has no successor in
$\widetilde{\Gamma}$ lying in $\mathcal X$, because $\mathcal X$ is
stable under predecessors in $\widetilde{\Gamma}_{\geqslant n}$. If
there is no arrow $y\to x$ in $\widetilde{\Gamma}$ such that
$y\in\mathcal X$, then $x$ has no predecessor in $\mathcal X$, by
minimality of $l(x)$. In such a case the full subquiver $\mathcal X'$
of $\widetilde{\Gamma}$ generated by $\mathcal X$ and $x$ has the same
arrows as those in $\mathcal X$ and  it is convex. Then $p$ trivially
extends to a well-behaved functor $p'\colon \ts k \mathcal
X'\to\ts{ind}\,\Gamma$ so that $(\mathcal X',p')\in\Sigma'$ and
$(\mathcal X,p)<(\mathcal X',p')$, a contradiction to the maximality
of $(\mathcal X,p)$. On the other hand, if there is an arrow $y\to x$
in $\widetilde{\Gamma}$ with $y\in\mathcal X$, then,  using dual
arguments to those used on the previous situation
(when $\mathcal Y$ was supposed to satisfy (a)), we similarly extend
$p$ to a well-behaved functor $p'\colon \ts k\mathcal X'\to
\ts{ind}\,\Gamma$ where $\mathcal X'$ is the (convex) full subquiver
of $\widetilde{\Gamma}$ generated by $\mathcal X$ and $x$. As in the
previous case, $(\mathcal X',p')\in\Sigma'$ and
$(\mathcal X,p)<(\mathcal X',p')$, which contradict the maximality of
($\mathcal X,p)$. Therefore, $p\colon \ts k\mathcal X\to \ts{ind}\,\Gamma$
is a well-behaved functor extending $q$, where $\mathcal X$ equals $\widetilde{\Gamma}_{\geqslant
  n}$ (and therefore  satisfies (c)).\\

Finally, assume that $\mathcal Y$ satisfies (c). Let $\Sigma''$ be the
set of those pairs $(\mathcal X,p)$ where $\mathcal X$ is a full
convex subquiver of $\widetilde{\Gamma}$ containing both
$\widetilde{\Gamma}_{\geqslant n}$ and $\mathcal Y$, and $\mathcal X$ is stable under
successors in $\widetilde{\Gamma}$ and $p\colon \ts k\mathcal
X\to\ts{ind}\,\Gamma$ is a well-behaved functor extending $q$. Then
$\Sigma''$ is non-empty for it contains $(\mathcal Y,q)$, and
$(\Sigma'', \leqslant)$ is totally inductive (with $\leqslant$ as
above). Let $(\mathcal X,p)$ be a maximal element in $\Sigma''$. We claim that
$\mathcal X=\widetilde{\Gamma}$. By absurd, assume that this is not
the case. Let $x\in\widetilde{\Gamma}$ be a vertex not in $\mathcal X$
and with $l(x)$ maximal for this property. This is possible because
$\widetilde{\Gamma}_{\geqslant n}\subseteq \mathcal X$. Since
$\mathcal X$ is stable under successors in $\widetilde{\Gamma}$, there
is no arrow $z\to x$ with $z\in\mathcal X$. Since $\widetilde{\Gamma}$
is connected, there exists an arrow $x\to y$ in
$\widetilde{\Gamma}$. The vertex $y$ then lies in $\mathcal X$ by
maximality of $l(x)$. Let $\mathcal X'$ be the full subquiver of
$\widetilde{\Gamma}$ generated by $\mathcal X$ and $x$. Therefore, the
arrows in $\mathcal X'$ are those in $\mathcal X$ together with those
in $\widetilde{\Gamma}$ starting in $x$ (which, by maximality of $l(x)$
have their endpoint in $\mathcal X$) and $\mathcal X'$ is
convex. Now, using the same arguments
as those used in the first situation (when we assumed that $\mathcal
Y$ satisfied (a)), we extend $p$ to a well-behaved functor $p'\colon
\ts k\mathcal X'\to\ts{ind}\,\Gamma$. We thus have $(\mathcal
X',p')\in\Sigma''$ and $(\mathcal X,p)<(\mathcal X',p')$, a
contradiction to the maximality of $(\mathcal X,p)$. This proves that
$\mathcal X=\widetilde{\Gamma}$ and finishes the proof of the
proposition.\end{proof}

We now present some practical situations where \ref{prop:well_behaved_general}
may be applied. 

\begin{defn}
\label{defn:sectional_family}
Let $X$ be an indecomposable module in $\Gamma$ and $r\geqslant 1$. A
\emph{sectional family of paths (starting in $X$ and of irreducible
  morphisms)} is a family
\begin{equation}
\xymatrix{
& X_{1,1}\ar[r]  &  \cdots \ar[r] & X_{1,l_1-1}
\ar[r]^{f_{1,l_1}} & X_{1,l_1}\\
X \ar[ru]^{f_{1,1}} \ar[r]^{f_{2,1}}  \ar[rd]_{f_{r,1}} & X_{2,1}
\ar[r] \ar@{.}[d]& \cdots \ar[r] & X_{2,l_2-1} 
\ar[r]^{f_{2,l_2}} & X_{2,l_2}\\
&X_{r,1} \ar[r]  & \cdots \ar[r] & X_{r,l_r-1}
\ar[r]_{f_{r,l_r}} & X_{r,l_r}
}\notag
\end{equation}
of $r$ paths starting in $X$ and of irreducible morphisms between
indecomposables, subject to the following conditions (where
$X=X_{i,0}$, for convenience):
\begin{enumerate}[(a)]
\item For every $M\in\Gamma$ and $l\geqslant 1$, let $I$ be the set of
  those indices $i\in\{1,\ldots,r\}$ such that $l_i\geqslant l$ and
  $f_{i,l}$ has domain $M$. Then the morphism $[f_{i,l}\ ;\ i\in
  I]\colon M\to \bigoplus\limits_{i\in I}X_{i,l}$ is irreducible.
\item For every $M\in\Gamma$ and $l\geqslant 1$, let $J$ be the set of
  those indices $i\in\{1,\ldots,r\}$ such that $l_i\geqslant l$ and
  $f_{i,l}$ has codomain $M$. Then the morphism $[f_{i,l}\ ;\ i\in
  I]^t\colon \bigoplus\limits_{i\in I}X_{i,l-1}\to M$ is irreducible.
\item There is no hook of the form $\cdot \xrightarrow{f_{i,j}}\cdot \xrightarrow{f_{i',j+1}}\cdot$.
\end{enumerate}
\end{defn}

\begin{rem}
\label{rem:sectional_family}
\begin{enumerate}[(1)]
\item The definition implies that each of 
the paths in the given family is sectional.
\item If $r=1$ the definition coincides with that of a sectional path.
\item If $l_i=1$ for every $i$, then the definition is equivalent to
  say that the morphism $[f_{1,1},\ldots,f_{r,1}]^t\colon X\to
  \bigoplus\limits_{i=1}^rX_{i,1}$ is irreducible.
\item Since $\ts k$ is an algebraically closed field, the fist two
  conditions together  are equivalent to the following single
  condition: For every $M,N\in\Gamma$ and $l\geqslant 1$, let $K$ be
  the set of indices $i\in\{1,\ldots,r\}$ such that $l_i\geqslant l$
  and $f_{i,l}$ is a morphism from $M$ to $N$, then the morphisms
  $f_{i,l}\colon M\to N$, $i\in K$, are linearly independent modulo
  $\ts{rad}^2(M,N)$.
\end{enumerate}
\end{rem}

\begin{prop}
  \label{prop:well_behaved_sectional}
Let $X$ be in $\Gamma$ and $x\in\pi^{-1}(X)$. Let
$\{X\xrightarrow{f_{i,1}}X_{i,1}\to\cdots\to
X_{i,l_i-1}\xrightarrow{f_{i,l_i}}X_{i,l_i}\}_{i=1,\ldots,r}$ be a
sectional family of paths starting in $X$ and of irreducible morphisms.
 Then there exist $r$ paths in $\widetilde{\Gamma}$
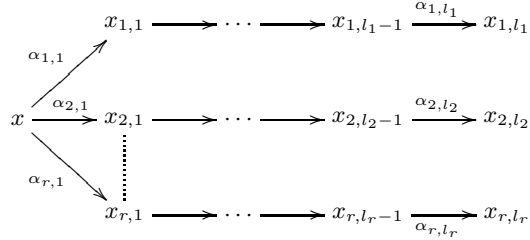
\begin{equation}
\xymatrix{
& x_{1,1}\ar[r]  &  \cdots \ar[r] & x_{1,l_1-1}
\ar[r]^{\alpha_{1,l_1}} & x_{1,l_1}\\
x \ar[ru]^{\alpha_{1,1}} \ar[r]^{\alpha_{2,1}}  \ar[rd]_{\alpha_{r,1}} & x_{2,1}
\ar[r] \ar@{.}[d]& \cdots \ar[r] & x_{2,l_2-1} 
\ar[r]^{\alpha_{2,l_2}} & x_{2,l_2}\\
&x_{r,1} \ar[r]  & \cdots \ar[r] & x_{r,l_r-1}
\ar[r]_{\alpha_{r,l_r}} & x_{r,l_r}
}\notag
\end{equation}
starting in $x$, such that $\pi x_{i,j}=X_{i,j}$, for every $i,j$, and
the arrows $\alpha_{i,j}$ are pairwise distinct. Moreover, for any
such data, there exists a well-behaved functor $F\colon
\ts k(\widetilde{\Gamma})\to\ts{ind}\,\Gamma$ such that
$F(\overline{\alpha_{i,j}})=f_{i,j}$, for every $i,j$.
\end{prop}
\begin{proof}
  We first construct the vertices $x_{i,j}$ and the arrows
  $\alpha_{i,j}$. For every $i\in\{1,\ldots,r\}$, the path $X\xrightarrow{f_{i,1}}X_{i,1}\to\cdots\to
X_{i,l_i-1}\xrightarrow{f_{i,l_i}}X_{i,l_i}$ of irreducible morphisms
defines a (non-unique) path $X\to X_{i,1}\to\cdots\to
X_{i,l_i-1}\to X_{i,l_i}$ in $\Gamma$. Since
$\pi\colon\widetilde{\Gamma}\to\Gamma$ is a covering of quivers, this
path in $\Gamma$ defines a path $x\to x_{i,1}\to\cdots\to
x_{i,l_i-1}\to x_{i,l_i}$ in $\widetilde{\Gamma}$ such that $\pi
x_{i,j}=X_{i,j}$. This defines all the vertices $x_{i,j}$. Let
$y,z\in\widetilde{\Gamma}$ be vertices and let $K$ be the set of
couples $(i,j)$, $i\in\{1,\ldots,r\}$ and $j\in\{1,\ldots,l_i\}$, such
that $y=x_{i,j-1}$ and $z=x_{i,j}$ (with the convention $x_{i,0}=x$). Note that if both $(i,j)$ and $(i',j')$ lie in $K$,
then $j=j'$ because $\widetilde{\Gamma}$ is with length and $x_{i,j}$
(or $x_{i',j'}$) is the endpoint of a path in $\widetilde{\Gamma}$
starting in $x$ and of length $j$ (or $j'$, respectively). By
definition of a sectional family of paths, the irreducible morphisms
$f_{i,j}\colon \pi y\to \pi z$, for $(i,j)\in K$, are linearly
independent modulo $\ts{rad}^2(\pi y,\pi z)$. Since
$\pi\colon\widetilde{\Gamma}\to\Gamma$ is a covering of quivers, there
is an injective map $(i,j)\mapsto \alpha_{i,j}$ from $K$ to the set of
arrows from $y$ to $z$ in $\widetilde{\Gamma}$. By proceeding this
construction for every vertices $y,z\in\widetilde{\Gamma}$, one
defines all the arrows $\alpha_{i,j}$, for $i\in\{1,\ldots,r\}$ and
$j\in\{1,\ldots,l_i\}$, which are pairwise distinct, by
construction.

Given the vertices $x_{i,j}$ and the arrows $\alpha_{i,j}$ as above,
we let $\mathcal Y$ be the full subquiver of $\widetilde{\Gamma}$
generated by all the $x_{i,j}$.  We need some properties on $\mathcal
Y$. Note that if there exists a path in $\widetilde{\Gamma}$ of the
form $x_{i,j}\to y_1\to\cdots\to y_s\to x_{i',j'}$, for some vertices
$y_1,\ldots,y_s\in\widetilde{\Gamma}$, then we have two parallel paths
in $\widetilde{\Gamma}$
\begin{align}
&  x\xrightarrow{\alpha_{i,1}}x_{i,1}\to\cdots\to
   x_{i,j-1}\xrightarrow{\alpha_{i,j}}x_{i,j} \to y_1\to\cdots \to y_s\to
   x_{i',j'}, \text{and}\notag\\
&x\xrightarrow{\alpha_{i',1}}x_{i',1}\to\cdots\to
   x_{i',j'-1}\xrightarrow{\alpha_{i',j'}}x_{i',j'}.\notag
\end{align}
Note that the image under $\pi$ of the second path is a path $X\to X_{i',1}\to\cdots\to
   X_{i',j'-1}\to X_{i',j'}$ in $\Gamma$ which is sectional
   (\ref{rem:sectional_family}, (1)). Since $\pi\colon
   \widetilde{\Gamma}\to\Gamma$ is a covering of translation quivers,
   this implies that the path $x\xrightarrow{\alpha_{i',1}}x_{i',1}\to\cdots\to
   x_{i',j'-1}\xrightarrow{\alpha_{i',j'}}x_{i',j'}$ is
   sectional. Applying \ref{prop:generic_covering} then shows
   that $j'=j+s+1$ and the sequence of vertices $(x,x_{i,1},\ldots,
   x_{i,j},y_1,\ldots,y_s,x_{i',j'})$ and
   $(x,x_{i',1},\ldots,x_{i',j'})$ coincide. From this, we deduce the
   following facts:
   \begin{enumerate}[(1)]
   \item $\mathcal Y$ is convex in $\widetilde{\Gamma}$.
   \item $\mathcal Y$ contains no path $y\to \cdot\to z$ with $z$
     non-projective and $y=\tau z$.
   \item If there is an arrow $y\to z$ in $\mathcal Y$, then there
     exist $i,j$ such that $y=x_{i,j-1}$ and $z=x_{i,j}$. Moreover,
     given $i',j'$, we have $y=x_{i',j'-1}$ if and only if $z=x_{i',j'}$.
   \end{enumerate}
We now define a well-behaved functor $q\colon \ts k\mathcal
Y\to\ts{ind}\,\Gamma$ such that $q(\alpha_{i,j})=f_{i,j}$ for every
$i,j$. Let $y,z\in\mathcal Y$ be vertices such that there exists at
least one arrow from $y$ to $z$. Then there is a path in
$\widetilde{\Gamma}$ from $x$ to $z$, say of length $n$. Since
$\widetilde{\Gamma}$ is with length and because of (3) above, we
deduce that if $z=x_{i,j}$ for some $i,j$, then $j=n$ and
$y=x_{i,j-1}$. We thus define $I_{y,z}$ to be the set of indices
$i\in\{1,\ldots,r\}$ such that $n\geqslant l_i$ and $y=x_{i,n-1},
z=x_{i,n}$. The set of arrows in $\widetilde{\Gamma}$ from $y$ to $z$
is therefore equal to $\{\alpha_{i,n}\ |\ i\in
I_{y,z}\}\cup\{\gamma_1,\ldots,\gamma_s\}$ where
$\gamma_1,\ldots,\gamma_s$ are pairwise distinct arrows, none of which
is equal to either of the arrows $\alpha_{i',j'}$,
$i'\in\{1,\ldots,r\}$ and $j'\in\{1,\ldots,l_{i'}\}$. Recall that the
irreducible morphisms $f_{i,n}\colon \pi y\to \pi z$, for $i\in
I_{y,z}$, are linearly independent modulo $\ts{rad}^2(\pi y,\pi z)$
(\ref{defn:sectional_family} and \ref{rem:sectional_family}). Since
$\pi$ induces a bijection from the set of arrows in
$\widetilde{\Gamma}$ from $y$ to $z$ to the set of arrows in $\Gamma$
from $\pi y$ to $\pi z$, we deduce that there exist irreducible
morphisms $g_1,\ldots,g_s\colon \pi y\to \pi z$ such that
$g_1,\ldots,g_s$ together with $f_{i,n}$, for $i\in I_{y,z}$, are
linearly independent modulo $\ts{rad}^2(\pi y,\pi z)$. We then set
$q(\alpha_{i,n})=f_{i,n}$, for every $i\in I_{y,z}$, and
$q(\gamma_j)=g_j$, for every $j=1,\ldots,s$. This defines $q$ on every
arrow from $y$ to $z$, for every vertices $y,z\in\mathcal Y$. Hence
the functor $q\colon \ts k \mathcal Y\to \ts{ind}\,\Gamma$. The
construction of $q$ and the above property (2) of $\mathcal Y$ show
that $q$ is well-behaved and $q(\alpha_{i,j})=f_{i,j}$, for every
$i,j$.

Finally, let $l$ be a length function on $\widetilde{\Gamma}$. Then
$\mathcal Y$ is a convex full subquiver of $\widetilde{\Gamma}$ such
that $l(y)\in\{l(x),l(x)+1,\ldots,l(x)+\underset{i=1,\ldots,r}{\ts{max}}
l_i\}$, for every vertex $y\in\mathcal Y$. Therefore,
\ref{prop:well_behaved_general}, (a), implies that there exists a
well-behaved functor $F\colon
\ts k(\widetilde{\Gamma})\to\ts{ind}\,\Gamma$ such that
$F(\overline{\alpha_{i,j}})=q(\alpha_{i,j})=f_{i,j}$, for every $i,j$.
\end{proof}

\begin{rem}
The proofs we gave for \ref{prop:well_behaved_general} and
\ref{prop:well_behaved_sectional} strongly rely on the fact that the
generic covering $\pi$ induces a bijection from the set of arrows in
$\widetilde{\Gamma}$ from $x$ to $y$ to the set of arrows in $\Gamma$
from $\pi x$ to $\pi y$, for every vertices $x,y\in \widetilde{\Gamma}$. In particular, these proofs are not likely to
be adapted to the situation where one replaces the generic covering
$\widetilde{\Gamma}$ of $\Gamma$ by the universal covering.
\end{rem}

The following result follows from
\ref{prop:well_behaved_sectional}. It will be particularly useful to us.
\begin{prop}
\label{lem:well_behaved}
  Let $X,X_1,\ldots,X_r$ lie on $\Gamma$ and
  $f=\begin{bmatrix}f_1,\ldots,f_r\end{bmatrix}^t\colon X\to
  X_1\oplus\ldots\oplus X_r$ be an irreducible morphism in $\ts{mod}\,
  A$. Let $x\in \pi^{-1}(X)$ and  $x\xrightarrow{\alpha_i}x_i$  be
  an arrow in $\widetilde{\Gamma}$ such that $\pi x_i=X_i$ for every
  $i\in\{1,\ldots,r\}$. Then there exists a well-behaved functor
  $F\colon \ts k(\widetilde{\Gamma})\to\ts{ind}\,\Gamma$ such that $F(\overline{\alpha_i})=f_i$ for every $i$.
\end{prop}
\begin{proof}
  It follows from \ref{rem:sectional_family}, (3), that the family
  of morphisms $\{f_1,\ldots,f_r\}$ is a sectional family of paths
  starting in $X$. The conclusion thus
  follows from \ref{prop:well_behaved_sectional}.
\end{proof}

We now study some properties of well-behaved functors which are
essential to our work. We begin with the following basic lemma.
\begin{lem}
\label{lem:covering}
Let $F\colon \ts k(\widetilde{\Gamma})\to\ts{ind}\,\Gamma$ be a well-behaved functor, $x,y$
vertices in $\widetilde{\Gamma}$ and $n\geqslant 0$. Then:
\begin{enumerate}[(a)]
\item $F$ maps a morphism in $\mathfrak{R}^n\ts k(\widetilde{\Gamma})(x,y)$ onto a morphism in
  $\ts{rad}^n(Fx,Fy)$.
\item Let $f\in\ts{rad}^{n+1}(Fx,Fy)$ and $\alpha_1\colon x\to
  x_1,\ldots, \alpha_r\colon x\to x_r$ be the arrows in
  $\widetilde{\Gamma}$ starting from $x$. Then there exist
  $h_i\in\ts{rad}^n(Fx_i,Fy)$, for every $i$, such that $f=\sum\limits_ih_iF(\overline{\alpha}_i)$.
\end{enumerate}
\end{lem}
\begin{proof} 
(a) follows from the fact that $F$ is well-behaved.

(b) We have a decomposition $f=\sum\limits_jg_jf_j$ where $j$ runs
through some index set, $f_j\in\ts{rad}(Fx,Y_j)$,
$g_j\in\ts{rad}^n(Y_j,Fy)$, $Y_j\in\ts{ind}\,A$, for every $j$. The
morphism $[f(\overline{\alpha}_1),\ldots,f(\overline{\alpha}_r)]^t\colon
Fx\to \bigoplus\limits_{i=1}^rF x_i$ is minimal left almost split
so every $f_j$ factors through it:
$f_j=\sum\limits_{i=1}^rf'_{j,i}F(\overline{\alpha}_i)$ with
$f'_{i,j}\in\ts{Hom}_A(Fx_i,Y_j)$. Setting
$h_i=\sum\limits_jg_jf'_{j,i}$ does the trick.\end{proof}

The following theorem states the main properties of well-behaved
functors we shall use.
 Part (b) of it was
first proved in \cite[\S 2]{R80a} in the case of the stable part of
the Auslander-Reiten quiver of a self-injective algebra of finite
representation type (see also \cite[3.1 Ex. (b)]{BG82} for the case of
the Auslander-Reiten quiver of an algebra of finite representation type).
\begin{Thm}
\label{prop:covering}
Let $F\colon \ts k(\widetilde{\Gamma})\to\ts{ind}\,\Gamma$ be a well-behaved functor, $x,y$
vertices in $\widetilde{\Gamma}$ and $n\geqslant 0$. Then:
\begin{enumerate}[(a)]
\item The two following maps induced by $F$ are bijective:
 \begin{equation}
   \begin{array}{rclc}
     \bigoplus\limits_{Fz=Fy}\mathfrak{R}^n\ts k(\widetilde{\Gamma})(x,z)
     /
\mathfrak{R}^{n+1}\ts k(\widetilde{\Gamma})(x,z)
     & \to &
     \ts{rad}^n(Fx,Fy)/\ts{rad}^{n+1}(Fx,Fy) \\
\\
     \bigoplus\limits_{Fz=Fy}\mathfrak{R}^n\ts k(\widetilde{\Gamma})(z,x)/
     \mathfrak{R}^{n+1}\ts k(\widetilde{\Gamma})(z,x)
     & \to &
     \ts{rad}^n(Fy,Fx)/\ts{rad}^{n+1}(Fy,Fx) & .
   \end{array}\notag
 \end{equation}
\item The two following maps induced by $F$ are injective:
  \begin{equation}
     \bigoplus\limits_{Fz=Fy}\ts k(\widetilde{\Gamma})(x,z) \to
     \ts{Hom}_A(Fx,Fy)
\ \ and\ \ 
     \bigoplus\limits_{Fz=Fy}\ts k(\widetilde{\Gamma})(z,x) \to 
     \ts{Hom}_A(Fy,Fx) .
\notag
  \end{equation}
\item $\Gamma$ is generalized standard if and only if $F$ is a covering
  functor, that is, the two maps of (b) are bijective (see \cite[3.1]{BG82}).
\end{enumerate}
\end{Thm}
\begin{proof} 
We prove the assertions concerning morphisms $Fx\to Fy$. Those
concerning $Fy\to Fx$ are proved using similar arguments. Let $\alpha_i\colon
x\to x_i$, $i=1,\ldots,r$, be the arrows in $\widetilde{\Gamma}$ starting from
$x$. So we have a minimal left almost split morphism in $\ts{mod}\, A$:
\begin{equation}
  Fx\xrightarrow{[F(\overline\alpha_1),\ldots,
      F(\overline\alpha_r)]^t}\bigoplus\limits_{i=1}^rFx_i\ .\notag
\end{equation}
(a) We denote by $F_n$ the map
$\bigoplus\limits_{Fz=Fy} \mathfrak{R}^n\ts k(\widetilde{\Gamma})(x,z)/
\mathfrak{R}^{n+1}\ts k(\widetilde{\Gamma})(x,z)\to
     \ts{rad}^n(Fx,Fy)/\ts{rad}^{n+1}(Fx,Fy)$. We prove that $F_n$ is
     surjective by induction on $n\geqslant 0$. So given a morphism
     $f\in\ts{rad}^n(Fx,Fy)$ we prove that there exists
     $(\phi_z)_z\in\bigoplus\limits_{Fz=Fy}\mathfrak{R}^n\ts k(\widetilde{\Gamma})(x,z)$ such that
     $f=\sum\limits_zF(\phi_z)\ \ts{mod}\,\ \ts{rad}^{n+1}$. We start
     with $n=0$. Let $f\in\ts{Hom}_A(Fx,Fy)$. If
$Fx\neq Fy$ then $f\in\ts{rad}(Fx,Fy)$. Otherwise, $f=\lambda 1_{Fx}\ \ts{mod}\,\
\ts{rad}$ with $\lambda\in \ts k$, that is, $f=F(\lambda 1_x)\ \ts{mod}\,\ \ts{rad}$ for some $\lambda\in
\ts k$. So $F_0$ is surjective. 
% Now we treat the case $n=1$. Let $f\in\ts{rad}(Fx,Fy)$. 
% So there is a factorization:
%   $f=\sum\limits_if_iF(\overline\alpha_i)$
% with $f_i\in\ts{Hom}_A(Fx_i,Fy)$ for every $i$, because of
% \ref{lem:covering}, (b). Let
% $i\in\{1,\ldots,r\}$. 
% If $f_i\in\ts{rad}(Fx_i,Fy)$, then
% $f_iF(\overline{\alpha}_i)\in\ts{rad}^2(Fx,Fy)$. Otherwise,
% $f_i=\lambda_iId_{Fx_i}\ \ts{mod}\ \ts{rad}$ with $\lambda_i\in
% k^*$. Thus
% $f=\sum\limits_{f_i\not\in\ts{rad}}F(\lambda_i\overline{\alpha}_i)\
% \ts{mod}\ \ts{rad}^2$.
%  This proves that $F_1$ is surjective.
Now let $n\geqslant 0$ and assume that $F_n$ is surjective. Let
$f\in\ts{rad}^{n+1}(Fx,Fy)$. 
Because of \ref{lem:covering}, (b), there is a decomposition
$f=\sum\limits_ih_iF(\overline{\alpha}_i)$ with
$h_i\in\ts{rad}^n(Fx_i,Fy)$. Moreover, 
$h_i=\sum\limits_zF(\phi_{i,z})\ \ts{mod}\ \ts{rad}^{n+1}$
with 
$(\phi_{i,z})_z\in\bigoplus\limits_{Fz=Fy}\mathfrak{R}^n\ts k(\widetilde{\Gamma})(x_i,z)$,
for every $i$, because $F_n$ is surjective. Therefore, 
  $f = \sum\limits_zF \left( \sum\limits_i\phi_{i,z}\overline\alpha_i
  \right)\ \ts{mod}\ \ts{rad}^{n+2}$
and
$\sum\limits_i\phi_{i,z}\overline{\alpha}_i\in\mathfrak{R}^{n+1}\ts k(\widetilde{\Gamma})(x,z)$,
for every $z\in\widetilde{\Gamma}$ such that $Fz=Fy$. So $F_{n+1}$ is
surjective. This proves that $F_n$ is surjective for every $n\geqslant 0$.\\

Now we prove that $F_n$ is injective for every  $n\geqslant
0$. 
Actually, we prove that the following assertion ($H_n$)
holds true:
\emph{``Let $(\phi_z)_z\in\bigoplus\limits_{Fz=Fy}\ts k(\widetilde{\Gamma})(x,z)$ be such that
 $\sum\limits_zF(\phi_z)\in\ts{rad}^n$, then
$\phi_z\in\mathfrak{R}^n\ts k(\widetilde{\Gamma})(x,z)$ for every $z$''.}
Clearly, this will prove the injectivity of all the $F_n$. We proceed
by induction on $n\geqslant 0$.
 Assume that $n=0$ and that
$\sum\limits_zF(\phi_z)\in\ts{rad}(Fx,Fy)$. If $Fx\neq Fy$ then $x\neq z$
for every $z$ such that $Fz=Fy$ and, therefore,
$\phi_z\in\mathfrak{R} \ts k(\widetilde{\Gamma})(x,z)$. If $Fx=Fy$ then $\phi_z\in\mathfrak{R} \ts k(\widetilde{\Gamma})(x,z)$ if
$x\neq z$ and there exists $\lambda\in \ts k$ such that $\phi_x=\lambda
1_x$. So $\lambda 1_{Fx}\in\ts{rad}(Fx,Fy)$, that is, $\lambda=0$. Thus,
$\phi_z\in\mathfrak{R} \ts k(\widetilde{\Gamma})(x,z)$ for every $z$. This proves
that ($H_0$) holds true. Now let $n\geqslant 0$, assume that ($H_n$)
holds true
and
let $(\phi_z)_z\in\bigoplus\limits_{Fz=Fy}\ts k(\widetilde{\Gamma})(x,z)$ be such that
$\sum\limits_zF(\phi_z)\in\ts{rad}^{n+2}$. So,
$\phi_z\in\mathfrak{R}^{n+1}\ts k(\widetilde{\Gamma})(x,z)$, for every
$z$, because ($H_n$) holds true. Also, 
there exists 
 $(\psi_z)\in\bigoplus\limits_{Fz=Fy}\mathfrak{R}^{n+2}\ts k(\widetilde{\Gamma})(x,z)$
such that
$\sum\limits_zF(\phi_z)=\sum\limits_zF(\psi_z)\ \ts{mod}\
\ts{rad}^{n+3}$, because $F_{n+2}$ is surjective and
$\sum\limits_zF(\phi_z)\in \ts{rad}^{n+2}(Fx,Fy)$.
Therefore, there exists $h_i\in \ts{rad}^{n+2}(Fx_i,Fy)$, for every $i$,
such that $\sum\limits_zF(\phi_z-\psi_z)=\sum\limits_i
  h_iF(\overline\alpha_i)$, because of \ref{lem:covering}, (b).
Since $\phi_z,\psi_z\in\mathfrak{R}\ts k(\widetilde{\Gamma})(x,z)$,  there is a decomposition
$\phi_z-\psi_z=\sum\limits_{i=1}^r\theta_{z,i}\overline\alpha_i$ with
   $\theta_{z,i}\in \ts k(\widetilde{\Gamma})(x_i,z)$ for every $i$.
We deduce that
\begin{equation}
\sum\limits_i \left(\sum\limits_zF(\theta_{z,i}) -
  h_i\right)F(\overline\alpha_i)=0\,\tag{$\star$}
\end{equation}
  Now if $x$ is injective, then
  $\sum\limits_z F(\theta_{z,i})-h_i=0$
 for every $i$. Since $h_i\in\ts{rad}^{n+2}(Fx_i,Fy)$, we deduce that
 $\sum\limits_zF(\theta_{z,i})\in\ts{rad}^{n+2}(Fx_i,Fy)\subseteq\ts{rad}^{n+1}(Fx_i,F_y)$
for every $i$. Because ($H_n$) holds true, we get
$\theta_{z,i}\in\mathfrak{R}^{n+1}\ts k(\widetilde{\Gamma})(x_i,z)$, for every
$i$, and, therefore
$\phi_z=\psi_z+\sum\limits_i\theta_{z,i}\overline{\alpha}_i\in\mathfrak{R}^{n+2}\ts k(\widetilde{\Gamma})$.
This proves that ($H_n$) holds true if $x$ is injective.
 Now
assume that $x$ is not injective. The mesh in $\widetilde{\Gamma}$ starting at $x$ is as follows:
\begin{equation}
  \xymatrix@R=1ex{
& x_1 \ar[rd]^{\beta_1} &\\
x \ar[ru]^{\alpha_1} \ar[rd]_{\alpha_r} & \vdots  & \tau^{-1}x\\ 
&x_r\ar[ru]_{\beta_r} & &.
}\notag
\end{equation}
Since $F$ is well-behaved, there is an almost split sequence in $\ts{mod}\,
A$:
\begin{equation}
  0\to Fx \xrightarrow{
     [ F(\overline\alpha_1) , \ldots ,F(\overline\alpha_r)]^t
} \bigoplus\limits_{i=1}^rFx_i \xrightarrow{
  [  F(\overline\beta_1) , \ldots , F(\overline\beta_r)]
} \tau_A^{-1}Fx\to 0\ .\notag
\end{equation}
From ($\star$), we deduce that there exists
$h\in\ts{Hom}_A(\tau_A^{-1}Fx,Fy)$ such that
$\sum\limits_zF(\theta_{z,i})-h_i=hF(\overline\beta_i)$, for every $i$. Since
    $F_0,\ldots,F_{n-1}$ are surjective, there exists
    $(\chi_z)_z\in\bigoplus\limits_{Fz=Fy}\ts k(\widetilde{\Gamma})(\tau^{-1}x,z)$ such that
    $h=\sum\limits_zF(\chi_z)\ \ts{mod}\,\ \ts{rad}^n$. Therefore, the
    following equality holds true for every $i$:
    \begin{equation}
      \sum\limits_z
      F(\theta_{z,i}) = \sum\limits_zF(\chi_z\overline\beta_i) +
      h_i\ \ts{mod}\ \ts{rad}^{n+1}\
      .\notag
    \end{equation}
Therefore, 
  $\sum\limits_z
  F(\theta_{z,i}-\chi_z\overline\beta_i)\in\ts{rad}^{n+1}(Fx_i,Fz)$,
  for every $i$, because $h_i\in\ts{rad}^{n+2}(Fx_i,Fy)$.
Hence,
$\theta_{z,i}-\chi_z\overline\beta_i\in\mathfrak{R}^{n+1}\ts k(\widetilde{\Gamma})(x_i,z)$, for every
$i,z$, because ($H_n$) holds true. This gives, for every $z$:
\begin{equation}
  \phi_z=\psi_z+\sum\limits_i
(\theta_{z,i}-\chi_z\overline\beta_i)
\overline\alpha_i\in\mathfrak{R}^{n+2}\ts k(\widetilde{\Gamma})(x,z)\ .\notag
\end{equation}
This proves that ($H_{n+1}$) holds true. Therefore, for every
$n\geqslant 0$ the map $F_n$ is injective and, therefore, bijective.\\

(b) Let $(\phi_z)_z\in\bigoplus\limits_{Fz=Fy}\ts k(\widetilde{\Gamma})(x,z)$ be such that
$\sum\limits_zF(\phi_z)=0$. In particular
$\sum\limits_zF(\phi_z)\in\ts{rad}^n(Fx,Fy)$ for every $n\geqslant
0$. Since $F_n$ is injective for every $n$, we deduce that
$\phi_z\in\mathfrak{R}^n \ts k(\widetilde{\Gamma})(x,z)$ for every $z$ and $n$. On the other
hand, given $z$ such that $Fz=Fy$, there exists $l\geqslant 0$ such
that all the paths from $x$ to $z$ in $\widetilde{\Gamma}$ are of length $l$, so that
$\mathfrak{R}^n\ts k(\widetilde{\Gamma})(x,z)=0$ for $n> l$. Therefore $\phi_z=0$ for every
$z$. This proves the injectivity of the first given map. The second
map is dealt with using dual arguments.\\

(c) Assume that $\Gamma$ is generalized standard and let
$f\in\ts{Hom}_A(Fx,Fy)$. So there
exists $n\geqslant 0$ such that $\ts{rad}^n(Fx,Fy)=0$. On the other hand,
the surjectivity of the maps $F_m$ ($m\geqslant 0$) shows that
$f=\sum\limits_zF(\phi_z)\ \ts{mod}\,\ \ts{rad}^n$ for
some $(\phi_z)_z\in\bigoplus\limits_{Fz=Fy}\ts k(\widetilde{\Gamma})(x,z)$. Therefore
$f=\sum\limits_zF(\phi_z)$. So the first given map in (c) is
surjective and so is the second one thanks to dual arguments. This and
(b) prove that $F$ is a covering functor.

Conversely, assume that $F$ is a covering functor and let $x,y\in\widetilde{\Gamma}$
be vertices. Therefore, there are only finitely many vertices $z\in\widetilde{\Gamma}$
such that $Fz=Fy$ and $\ts k(\widetilde{\Gamma})(x,z)\neq 0$ because $\ts{Hom}_A(Fx,Fy)$
is finite dimensional. This and the fact that $\widetilde{\Gamma}$ is with length
imply that  there exists $n\geqslant 0$ such that $\mathfrak{R}^n\ts k(\widetilde{\Gamma})(x,z)=0$
for every $z$ such that $Fz=Fy$. The injectivity of $F_n$ then
implies that $\ts{rad}^n(Fx,Fy)=0$. So $\Gamma$ is generalized standard.
\end{proof}

\begin{rem}
  It is not difficult to check that the proofs of
  \ref{lem:well_behaved} and Theorem~\ref{prop:covering}
  still work if $\Gamma$ is an Auslander-Reiten component of $\mathcal{T}$
  (instead of $\ts{mod}\,A$) where $\mathcal{T}$ is a triangulated
  Krull-Schmidt category over $\ts k$ with finite dimensional
  $\ts{Hom}$ spaces and Auslander-Reiten triangles.
\end{rem}

\section{Degrees of irreducible morphisms}
\label{sec:s2}

In this section we prove some characterizations for the left (or
right) degree of an irreducible morphism to be finite. These shall be
used later for the proof of our main result. Each statement has its
dual counterpart which will be omitted. The following proposition was
first proved in \cite{CPT04} for generalized standard convex
Auslander-Reiten components of an artin algebra. In a
weaker form it was also proved in \cite{C08} for standard
Auslander-Reiten components. We thank Shiping Liu for pointing out
that the arguments used to prove the first statement can be adapted
to prove the second statement. Note that the two statements are not dual
to each other.
\begin{prop}
\label{prop:degree}
  Let $f\colon X\to Y$ be an
  irreducible morphism with $X$ indecomposable, $\Gamma$ be the
  Auslander-Reiten component of $A$ containing $X$ and
  $n\in\mathbb{N}$. 
\begin{enumerate}[(a)]
\item If $d_l(f)=n$, then  there
  exist $Z\in\Gamma$ and $h\in\ts{rad}^n(Z,X)\backslash\ts{rad}^{n+1}(Z,X)$ such
  that $fh=0$.
\item [(b)] If $d_r(f)=n$, then  there
  exist $Z\in\Gamma$ and $h\in\ts{rad}^n(Y,Z)\backslash\ts{rad}^{n+1}(Y,Z)$ such
  that $hf=0$.
\end{enumerate}
\end{prop}
\begin{proof} 
We write $Y=X_1\oplus\cdots\oplus X_r$ with $X_1,\ldots,X_r\in\Gamma$
and $f=[f_1,\ldots,f_r]^t$ with $f_i\colon X\to X_i$. Let
$\pi\colon\widetilde{\Gamma}\to\Gamma$ be the generic
covering. Because $f$ is irreducible, $\widetilde{\Gamma}$ contains a
subquiver of the form
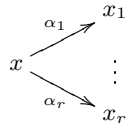
\begin{equation}
  \xymatrix@R=1ex{
&x_1\\
x\ar[ru]^{\alpha_1}\ar[rd]_{\alpha_r}&\vdots \\
&x_r &
}\notag
\end{equation}
such that $\pi x_i=X_i$ for every $i$.
Let $F\colon \ts k(\widetilde{\Gamma})\to\ts{ind}\, \Gamma$ be a
well-behaved functor such
that $F(\overline\alpha_i)=f_i$ for every $i$
(\ref{lem:well_behaved}).

(a)
If $d_l(f)=n$, then there exists $Z\in\Gamma$ and
$g\in\ts{rad}^n(Z,X)\backslash\ts{rad}^{n+1}(Z,X)$ such that
$fg\in\ts{rad}^{n+2}(Z,Y)$, that is $f_ig\in\ts{rad}^{n+2}(Z,X_i)$ for every $i$. 
Because of Theorem~\ref{prop:covering}, there exists $(\phi_z)_z\in\bigoplus\limits_{Fz=Z}\mathfrak{R}^n\ts k(\widetilde{\Gamma})(z,x)$ such that
$g=\sum\limits_zF(\phi_z)\
  \ts{mod}\,\ \ts{rad}^{n+1}(Z,X)$ and $\phi_{z_0}\not\in\mathfrak{R}^{n+1}\ts k(\widetilde{\Gamma})(z_0,x)$ for some $z_0$.
Therefore
$f_ig=\sum\limits_zF(\overline\alpha_i\phi_z)\
\ts{mod}\,\ \ts{rad}^{n+2}(Z,X_i)$ for every $i$.
Since $f_ig\in\ts{rad}^{n+2}(Z,X_i)$ we infer, using
 Theorem~\ref{prop:covering}, that
 $\overline{\alpha_i}\phi_z\in\mathfrak{R}^{n+2}\ts k(\widetilde{\Gamma})(z,x)$
 for every $z$ and every $i$. On the other hand,
 $\phi_{z_0}\not\in\mathfrak{R}^{n+1}\ts k(\widetilde{\Gamma})(z_0,x)$ implies that any path
 in $\widetilde{\Gamma}$ from $z_0$ to $x$ has length at most
 $n$. Hence, any path from $z_0$ to $x_i$ has length at most $n+1$ for
 every $i$. Thus $\overline{\alpha}_i\phi_{z_0}=0$ for every $i$. We
 then set $h=F(\phi_{z_0})$. Then
 $fh=\sum\limits_iF(\overline{\alpha}_i\phi_{z_0})=0$ and
   $h\in\ts{rad}^n(Z,X)\backslash\ts{rad}^{n+1}(Z,X)$, because
   $\phi_{z_0}\in\mathfrak{R}^n\backslash\mathfrak{R}^{n+1}$ and
   because of Theorem~\ref{prop:covering}.\\

(b) Now assume that $d_r(f)=n$. There exists $Z\in\Gamma$ and
$g\in\ts{rad}^n(Y,Z)\backslash\ts{rad}^{n+1}(Y,Z)$ such that
$gf\in\ts{rad}^{n+2}(X,Z)$. We write $g=[g_1,\ldots,g_r]$ with
$g_i\colon X_i\to Z$. Hence, $g_i\in\ts{rad}^n(X_i,Z)$; there exists
$i_0\in\{1,\ldots,r\}$ such that $g_{i_0}\not\in\ts{rad}^{n+1}(X_{i_0},Z)$;
  and $\sum\limits_ig_if_i\in\ts{rad}^{n+2}(X,Z)$.

For every $i$, there exists
$(\phi_{i,z})_z\in\bigoplus\limits_{Fz=Z}\mathfrak{R}^n\ts k(\widetilde{\Gamma})(x_i,z)$
such that $g_i=\sum\limits_zF(\phi_{i,z})\ \ts{mod}\ \ts{rad}^{n+1}$;
and, also, there exists $z_0$ such that $Fz_0=Z$ and
$\phi_{i_0,z_0}\not\in\mathfrak{R}^{n+1}\ts k(\widetilde{\Gamma})(x_{i_0}, z_0)$,
because of Theorem~\ref{prop:covering} and the above properties of the
$g_i$. In particular, the paths in $\widetilde{\Gamma}$ from $x_{i_0}$
to $z_0$ all have length at most $n$, and, therefore, the paths from
$x$ to $z_0$ all have length at most $n+1$.

On the other hand,
$\sum\limits_ig_if_i=\sum\limits_zF\left(\sum\limits_i\phi_{i,z}\overline{\alpha}_i\right)\
\ts{mod}\ \ts{rad}^{n+2}$ lies in $\ts{rad}^{n+2}(X,Z)$. Hence,
$\sum\limits_i\phi_{i,z}\overline{\alpha}_i\in\mathfrak{R}^{n+2}
\ts k(\widetilde{\Gamma})(x,z)$ for every $z$, because of
Theorem~\ref{prop:covering}. This and the above property on the length of the
paths in $\widetilde{\Gamma}$ from $x$ to $z_0$ imply that
$\sum\limits_i\phi_{i,z_0}\overline{\alpha}_i=0$ We then set
$h_i=F(\phi_{i,z_0})\colon X_i\to Z$ and $h=[h_1,\ldots,h_r]\colon
Y\to Z$. Then
$hf=F\left(\sum\limits_i\phi_{i,z_0}\overline{\alpha}_i\right)=0$,
$h\in\ts{rad}^n(Y,Z)$ because
$\phi_{i,z_0}\in\mathfrak{R}^n\ts k(\widetilde{\Gamma})(x_i,z_0)$  for
every $i$, and $h\not\in\ts{rad}^{n+1}(Y,Z)$ because
$\phi_{i_0,z_0}\not\in\mathfrak{R}^{n+1}
\ts k(\widetilde{\Gamma})(x_{i_0},z_0)$ (see Theorem~\ref{prop:covering}).
\end{proof}

 \begin{rem}
\label{rem:domain_explode}
   Keep the notations of \ref{prop:degree}. 
\begin{enumerate}[(a)]
\item
If $d_l(f)=n$, then, by definition, there exist $Z\in\Gamma$
   and $g\in\ts{rad}^n(Z,X)\backslash\ts{rad}^{n+1}(Z,X)$ such that
   $fg\in\ts{rad}^{n+2}(Z,X_1\oplus\cdots\oplus X_r)$. The proof of
   \ref{prop:degree} shows that there exists
   $h\in\ts{rad}^n(Z,X)\backslash\ts{rad}^{n+1}(Z,X)$ such that $fh=0$
   (that is, the domain of $h$ is equal to the domain of $g$).
Of course, the same remark holds true if $d_r(f)=n$.
\item It is still an open question to know whether the morphism
  $h$ in \ref{prop:degree} can be chosen to be a composition of
  irreducible morphisms (instead of a sum of compositions of
  such). Recall that this is indeed the case if $\alpha(\Gamma)\leqslant 2$
  (\cite{CPT04}).
\end{enumerate}
 \end{rem}

Now we derive some consequences of \ref{prop:degree}. 
The following corollary follows directly from
\ref{prop:degree}. We omit its proof. Note that it was proved in
\cite{C08} for irreducible morphisms between indecomposable modules
lying in a standard component.
\begin{cor}
  \label{cor:degree_sink}
Let $f\colon X\to Y$ be an
irreducible morphism in $\ts{mod}\, A$ with $X$ indecomposable. If
$d_l(f)$ is finite,
then
 $f$ is not mono and $d_r(f)=\infty$. In particular, every
 minimal left almost split morphism in $\ts{mod}\, A$
has infinite left degree.
\end{cor}

The following proposition compares $d_l(f)$ and $d_l(g)$
when there is an almost split sequence of the form $0\to
\tau_AY\xrightarrow{[g,g']^t} X'\oplus X\xrightarrow{[f',f]}Y\to
0$. Recall that it was proved in \cite[1.2]{L92} that $d_l(f)<\infty$
implies $d_l(g)\leqslant d_l(f)-1$ (in the more general setting of
artin algebras). Note that the following result was proved in
\cite{C08} in the case where the indecomposable module $Y$ lies in a
standard component.
\begin{prop}
\label{prop:shiping}
  Let $f\colon X\to Y$ be an irreducible morphism with $Y$
  indecomposable  and non-projective. Assume that the almost split sequence in
  $\ts{mod}\, A$:
  \begin{equation}
\xymatrix@R=1ex{
&&X'\ar[rd]^{f'}\\
0 \ar[r] & \tau_AY \ar[ru]^g \ar[rd]_{g'} &&Y \ar[r] & 0\\
&&X\ar[ru]_f&&
}\notag
\end{equation}
is such that $X'\neq 0$. Then $d_l(f)<\infty$ if and only if
$d_l(g)<\infty$. In such a case, $d_l(f)=n$ if and only if $d_l(g)=n-1$.
\end{prop}
\begin{proof} It was proved in \cite[1.2]{L92} that if
$d_l(f)<\infty$, then  $d_l(g)\leqslant
d_l(f)-1$. Converly, assume that $d_l(g)=m<\infty$. Then there exists
$Z\in\ts{ind}\,A$ and
$h\in\ts{rad}^m(Z,\tau_AY)\backslash\ts{rad}^{m+1}(Z,\tau_AY)$ such
that $gh=0$, because of \ref{prop:degree}. Consider the morphism
$g'h\in\ts{rad}^{m+1}(Z,X)$. The morphism $[g,g']^t$ is minimal left
almost split so it has infinite left degree, because of
\ref{cor:degree_sink}. Since $[g,g']^th=[0,g'h]$ we deduce that
$g'h\not\in\ts{rad}^{m+2}(Z,X)$. On the other hand,
$fg'h=(fg'+f'g)h=0$. This proves that 
$d_l(f)\leqslant m+1=d_l(g)+1$.\end{proof}

The following proposition is a key-step towards Theorem~\ref{thm1}.
\begin{prop}
\label{prop:kernel}
 Let $f\colon X\to Y$ be
an irreducible morphism   in $\ts{mod}\, A$ with
$X$ indecomposable, $\Gamma$ be the Auslander-Reiten component of $A$
containing $X$ and $n\geqslant 1$ be an integer. 
The two following conditions are equivalent:
\begin{enumerate}[(a)]
\item $d_l(f)=n$.
\item $f$ is not mono and the morphism $\ts{ker}(f)\colon
  \ts{Ker}(f)\hookrightarrow X$
  lies in
  $\ts{rad}^n(\ts{Ker}(f),X)\backslash\ts{rad}^{n+1}(\ts{Ker}(f),X)$.
\end{enumerate}
These conditions imply the following one:
\begin{enumerate}[(a)]
\setcounter{enumi}{2}
\item $f$ is not mono and $\ts{Ker}(f)\in\Gamma$.
\end{enumerate}
If $\Gamma$ is generalized standard, then the three conditions are equivalent.
\end{prop}
\begin{proof} 
If $d_l(f)=n<\infty$, then \ref{prop:degree}
implies that there exists $n\geqslant 0$, $Z\in\Gamma$ and
$h\in\ts{rad}^n(Z,X)\backslash\ts{rad}^{n+1}(Z,X)$ such that
$fh=0$. In particular, $f$ is not mono (and, therefore,
$\ts{Ker}(f)$ is indecomposable, because $f$ is
irreducible). Therefore, we have a factorization
\begin{equation}
  \xymatrix@R=2ex{
& Z \ar[d]^h \ar[ld]_{\exists}\\
\ts{Ker}(f)\ar@{^(->}[r] & X \ar[r]^(0.3)f & Y
}\notag
\end{equation}
which implies that $\ts{ker}(f)\not\in\ts{rad}^{n+1}(\ts{Ker}(f),X)$ and,
therefore, $\ts{Ker}(f)\in\Gamma$.
Let $i$ be such that $\ts{ker}(f)\in\ts{rad}^i(\ts{Ker}(f),X)$. So
$i\leqslant n$ and, since
$f\ts{ker}(f)=0$, we have $d_l(f)\leqslant i$. Thus, $i=n$ and
$\ts{ker}(f)\in\ts{rad}^n(\ts{Ker}(f),X)\backslash\ts{rad}^{n+1}(\ts{Ker}(f),X)$.
This proves that (a) implies (b) and (c).

If
$f$ is not mono and
$\ts{ker}(f)\in\ts{rad}^n(\ts{Ker}(f),X)\backslash\ts{rad}^{n+1}(\ts{Ker}(f),X)$
 then $\ts{Ker}(f)\in\Gamma$. From the equality
$f\ts{ker}(f)=0$ we deduce that $d_l(f)\leqslant n<\infty$. Since (a)
implies (b) we deduce that $d_l(f)=n$. This proves that (b) implies
(a) and (c).

Finally,  if $\Gamma$ is
generalized standard and  $\ts{Ker}(f)\in\Gamma$, then the inclusion
 $\ts{Ker}(f)\hookrightarrow X$ lies on
$\ts{rad}^n(\ts{Ker}(f),X)\backslash\ts{rad}^{n+1}(\ts{Ker}(f),X)$ for
some $n\geqslant 1$
because  $\ts{rad}^{\infty}(\ts{Ker}(f),X)=0$. Thus, (b) and,
therefore, (a) holds true.
\end{proof}

  Keep the notations of \ref{prop:kernel} and of its proof and assume
  that $d_l(f)=n$. Let
   $g\colon
  Z\to \ts{Ker}(f)$ a morphism such that $\ts{ker}(f) g=h$. 
Both morphisms $\ts{ker}(f)$ and $h$ lie in
$\ts{rad}^n\backslash\ts{rad}^{n+1}$ so that
$g\not\in\ts{rad}(Z,\ts{Ker}(f))$. Since both $Z$ and $\ts{Ker}(f)$
are indecomposable, we deduce that $g\colon Z\to \ts{Ker}(f)$ is an
isomorphism. In other words, we have the following

\begin{cor}
\label{rem:kernel_explode}
Let $f\colon X\to Y$ be an irreducible morphism with $X$
indecomposable. If $d_l(f)=n$ and if there exists $Z\in\ts{ind}\,A$
and $h\in\ts{rad}^n(Z,X)\backslash\ts{rad}^{n+1}(Z,X)$ such that
$fh=0$, then $h=\ts{ker}(f)$.
\end{cor}
\begin{proof} This follows from the arguments given before
the lemma.\end{proof}

% \begin{rem}
% \label{rem:cokernel}
%   Using similar arguments as those used in the proof of
%   \ref{prop:kernel}, one gets the following symmetric (but not dual)
%   result. Let $f\colon X\to Y$ be an irreducible morphism with
%   $X\in\ts{ind}\,A$, let $\Gamma$ be the Auslander-Reiten 
%  component of $A$
% containing $X$ and $n\geqslant 1$ be an integer. 
% The two following conditions are equivalent:
% \begin{enumerate}[(a)]
% \item $d_r(f)=n$.
% \item $f$ is not epi and the morphism $\ts{cokker}(f)\colon
%   Y\twoheadrightarrow \ts{Coker}(f)$
%   lies in
%   $\ts{rad}^n\backslash\ts{rad}^{n+1}$.
% \end{enumerate}
% These conditions imply the following one:
% \begin{enumerate}[(a)]
% \setcounter{enumi}{2}
% \item $f$ is not epi and $\ts{Coker}(f)\in\Gamma$.
% \end{enumerate}
% If $\Gamma$ is generalized standard, then the three conditions are
% equivalent.

% Moreover, as in \ref{rem:kernel_explode}, if $d_r(f)=n$ and if
% $Z\in\ts{ind}\,A$ and
% $h\in\ts{rad}^n(Y,Z)\backslash\ts{rad}^{n+1}(Y,Z)$ are such that
% $hf=0$, then $h=\ts{coker}(f)$.
% \end{rem}

Using \ref{prop:kernel} we can prove the following result.
\begin{cor}
  \label{cor:kernel_invariant}
  Let $f,f'\colon X\to Y$ be irreducible
  morphisms in $\ts{mod}\, A$ with $X$
  indecomposable. Then:
\begin{enumerate}[(a)]
\item 
If $f$ has finite left degree then $d_l(f)=d_l(f')$ and
  $\ts{Ker}(f)\simeq \ts{Ker}(f')$.
\item If $f$ has finite right degree then $d_r(f)=d_r(f')$.
% and
%   $\ts{Coker}(f)\simeq \ts{Coker}(f')$.
 \end{enumerate}
\end{cor}
\begin{proof} (a) Write $Y=X_1\oplus\cdots\oplus X_r$ with
$X_1,\ldots,X_r$ indecomposable.
Let $f_i,f_i'\colon X\to X_i$
($i\in\{1,\ldots,r\}$) be the
morphisms such that $f=\begin{bmatrix}f_1&\ldots&f_r\end{bmatrix}^t$ and
$f'=\begin{bmatrix}f_1' & \ldots & f_r'\end{bmatrix}^t$. By
\cite[Lem. 1.3]{L96}, the irreducible morphisms $f_1,\ldots,f_r$ all
have finite left degree. By \cite[Lem. 1.7]{L92} we deduce that for
every $i$ there exist a scalar $\lambda_i\in \ts k^*$ and a morphism
$r_i\in\ts{rad}^2(X,X_i)$ such that $f_i'=\lambda_if_i+r_i$. This clearly implies
 that $d_l(f)=d_l(f')$. Let $n=d_l(f)$ and $\iota\colon
\ts{Ker}(f)\hookrightarrow X$ be the inclusion. By \ref{prop:kernel}
we know that
$\iota\in\ts{rad}^n(\ts{Ker}(f),X)\backslash\ts{rad}^{n+1}(\ts{Ker}(f),X)$. On the
other hand, we have $f_i'\iota=r_i\iota\in\ts{rad}^{n+2}(\ts{Ker}(f),X_i)$ for
every $i$, that is,
$f'\iota\in\ts{rad}^{n+2}(\ts{Ker}(f),Y)$. By
\ref{prop:degree} and \ref{rem:domain_explode} we infer that there
exists $h\in\ts{rad}^n(\ts{Ker}(f),X)\backslash\ts{rad}^{n+1}(\ts{Ker}(f),X)$
such that $f' h=0$. Finally, \ref{rem:kernel_explode} implies
that $\ts{Ker}(f)\simeq\ts{Ker}(f')$.
\\

(b) If $X$ is injective, then there exist $U\in\ts{mod}\,A$ and
morphisms $u,u'\colon X\to U$ such that both $[f,u]$ and $[f',u']$ are
minimal right almost split morphisms $X\to Y\oplus U$. The dual
version of \ref{cor:degree_sink} implies that both $[f,u]$ and
$[f',u']$ have infinite right degree and, therefore, so do $f$ and
$f'$. Therefore, $X$ is not injective and there are minimal almost
split sequences in $\ts{mod}\,A$
\begin{equation}
  \xymatrix@R=5pt{
&Y\ar[rd]^g & && && & Y \ar[rd]^{g'} &\\
X \ar[ru]^f \ar[rd]_h&& \tau_A^{-1}X && \text{and} && X \ar[ru]^{f'}
\ar[rd]_{h'}&& \tau_A^{-1}X\\
&Y' \ar[ru]_{i} & && && & Y' \ar[ru]_{i'}&&.
}\notag
\end{equation}
The dual version of \ref{prop:shiping} applied to the first sequence yields 
$d_r(i)=d_r(f)-1$. Then, the dual version of (a) applied to
$i,i'\colon Y'\to \tau_A^{-1}X$ gives $d_r(i)=d_r(i')$. Finally, the
dual version of \ref{prop:shiping} applied to the second sequence
above yields $d_r(i')=d_r(f')-1$. Thus $d_r(f')=d_r(f)$.
\end{proof}

The following example shows that \ref{cor:kernel_invariant} does not
necessarily hold if one drops the finiteness condition on the left
degree.
\begin{ex}
  Let $A$ be the path algebra of the following quiver of type $\widetilde
  A_2$:
  \begin{equation}
    \xymatrix@R=1ex{
&2\ar[rd] \\
1\ar[ru] \ar[rr] & & 3 &.
}\notag
  \end{equation}
Given a vertex $x$, we write $I_x$ for the corresponding
indecomposable injective $A$-module. So the canonical quotient
$f\colon I_3\twoheadrightarrow  I_1$ is an irreducible morphism of
infinite left degree (see \cite[Cor. 4.10]{CPT04} for instance). Then
$\ts{Ker}(f)$ is as follows:
\begin{equation}
  \xymatrix@R=1ex{
&\ts k\ar[rd]^{\ts{Id}} \\
\ts k\ar[ru]^{\ts{Id}} \ar[rr]_0 & & \ts k &.
}\notag
\end{equation}
On the other hand, let $\mu\in\ts{rad}^2(I_3,I_1)$ be the composition
$I_3\twoheadrightarrow I_2\twoheadrightarrow I_1$ of the two canonical
quotients. Then $f'=f+\mu\colon I_3\twoheadrightarrow I_1$ is also irreducible and its
kernel is as follows:
\begin{equation}
  \xymatrix@R=1ex{
&\ts k\ar[rd]^{\ts{Id}} \\
\ts k\ar[ru]^{\ts{Id}} \ar[rr]_{\ts{Id}} & & \ts k &.
}\notag
\end{equation}
Clearly, $\ts{Ker}(f)$ and $\ts{Ker}(f')$ lie in distinct homogeneous tubes
and are therefore non-isomorphic.
\end{ex}

The result below follows from \ref{prop:kernel}. It was first
proved for standard algebras in \cite{C08}.
\begin{cor}
\label{cor:degree_finite_type}
  Let $A$ be of finite representation type and $f\colon X\to Y$  an
  irreducible morphism with $X$ or $Y$ indecomposable. Then the following
  conditions are equivalent:
  \begin{enumerate}[(a)]
  \item $d_l(f)<\infty$.
\item $d_r(f)=\infty$.
\item $f$ is an epimorphism.
  \end{enumerate}
\end{cor}
\begin{proof} 
Assume first that $X$ is indecomposable.
If $d_l(f)<\infty$, then $f$ is not
mono (and, therefore, it is epi, because it is
irreducible) and $d_r(f)=\infty$, because of the dual version of
\ref{cor:degree_sink}. So
(a) implies (b) and (c).

If $d_r(f)<\infty$, then $f$ is not epi, because of
\ref{prop:degree}. Therefore, (c) implies (b).
Conversely, if $f$ is not epi, then
$\ts{Coker}(f)\in\Gamma(\ts{mod}\,A)$ and there exists an integer
$n\geqslant 1$ such that $\ts{coker(f)}\colon
Y\twoheadrightarrow\ts{Coker}(f)$ lies in
$\ts{rad}^n(Y,\ts{Coker}(f))\backslash\ts{rad}^{n+1}(Y,\ts{Coker}(f))$
(indeed, $\ts{rad}^{\infty}=0$ because $A$ is
representation-finite). Thus, (b) implies (c) and these two conditions
are therefore equivalent.

Finally, assume that $f$ is epi. Then $f$ is not mono,
$\ts{Ker}(f)\in\Gamma(\ts{mod}\,A)$ and, as above, $\ts{ker}(f)\colon
\ts{Ker}(f)\hookrightarrow X$ lies in
$\ts{rad}^n(\ts{Ker}(f),X)\backslash\ts{rad}^{n+1}(\ts{Ker}(f),X)$. In
particular, $d_l(f)<\infty$. Thus (c) implies (a) and, therefore, the
three conditions are equivalent if $X$ is indecomposable.

If $Y$ is indecomposable, then, using dual arguments, one proves that the following
conditions are equivalent: $d_l(f)=\infty$; $d_r(f)<\infty$; and $f$
is mono. Since an irreducible morphism is either mono or
epi, this proves that (a), (b) and (c) are equivalent.\end{proof}

We end this section with another application of \ref{prop:degree}: The
description of the irreducible morphisms with indecomposable domain or
indecomposable codomain and with (left or right) degree equal to
$2$. Again, each statement has its dual counterpart which is
omitted. We thus restrict our study to irreducible morphisms with
indecomposable domain. We start with a characterization of the
equality $d_r(f)=2$.

\begin{cor}
  \label{cor:degree_right_two}
Let $f\colon X\to Y$ be an irreducible morphism with $X$
indecomposable. The following conditions are equivalent:
\begin{enumerate}[(a)]
\item $d_r(f)=2$.
\item $X$ is not injective and there exists an almost split
  sequence $0\to X\xrightarrow{[f,f']^t}Y\oplus
  Y'\xrightarrow{[g,g']} \tau_A^{-1}X\to 0$ with $Y'$ indecomposable
  non-injective fitting into an almost split sequence $0\to
  Y'\xrightarrow{g'}\tau_A^{-1}X\xrightarrow{\delta} \tau_A^{-1}Y'\to
  0$. In other words, there is a configuration of almost
  split sequences in $\ts{mod}\,A$
  \begin{equation}
    \xymatrix@!C@=3ex{
& Y\ar[rd]^g \\
X \ar[ru]^f \ar[rd]_{f'} \ar@{.}[rr] && \tau_A^{-1}X
\ar[rd]^{\delta}\\
& Y'  \ar[ru]_{g'}\ar@{.}[rr] && \tau_A^{-1}Y'&.
}\notag
  \end{equation}
\end{enumerate}
\end{cor}
\begin{proof} Assume that (a) holds true. Then $X$ is not
injective (\cite[1.3]{L92}) and $f$ is not a minimal left almost split
monomorphism
(\cite[1.12]{L92}). So there is an almost split sequence 
$0\to X\xrightarrow{[f,f']^t}Y\oplus
  Y'\xrightarrow{[g,g']} \tau_A^{-1}X\to 0$ with $Y'\in\ts{mod}\,A$
  non-zero. On the other hand, there exists $M\in\ts{ind}\,A$ and
  $h\in\ts{rad}^2(Y,M)\backslash\ts{rad}^3(Y,M)$ such that $hf=0$,
  because $d_r(f)=2$ and because of the dual version of \ref{prop:degree}. Since
  $[h,0][f,f']^t=0$, there exists $h'\colon\tau_A^{-1}X\to M$ such
  that $h'g=h$ and $h'g'=0$. Clearly,
  $h'\in\ts{rad}(\tau_A^{-1}X,M)\backslash\ts{rad}^2(\tau_A^{-1}X,M)$,
  because $g$ is irreducible and
  $h=h'g\in\ts{rad}^2(Y,M)\backslash\ts{rad}^3(Y,M)$. Hence,
  $d_r(g')=1$. Using \cite[1.12]{L92}, we deduce that $Y'$ is
  indecomposable not
  injective and $g'$ is minimal left almost split. This proves that
  (b) holds true.

Conversely, assume that (b) holds true. We prove that so does (a). In
particular, $X$ is non-injective and $f$ is not minimal left almost
split. Using \cite[1.12]{L92}, we infer that $d_r(f)\geqslant
2$. Consider the morphism $\delta
g\in\ts{rad}^2(Y,\tau_A^{-1}Y')$. Let $Z$ be an indecomposable summand
of $Y$ and let $Z\to \tau_A^{-1}X$ be the composition of $g\colon
Y\to\tau_A^{-1}X$ with the section $Z\to Y$. Then $Z\not\simeq Y'$,
because the almost split sequence starting at $Y'$  has its middle term
indecomposable. Therefore, $Z\to
\tau_A^{-1}X\xrightarrow{\delta}\tau_A^{-1}Y'$ is a sectional path of
irreducible morphisms so that its composite lies in
$\ts{rad}^2(Z,\tau_A^{-1}Y')\backslash\ts{rad}^3(Z,\tau_A^{-1}Y')$
(\cite{IT84}). Thus, $\delta g\in
\ts{rad}^2(Y,\tau_A^{-1}Y')\backslash\ts{rad}^3(Y,\tau_A^{-1}Y')$ and
$\delta gf=-\delta g'f'=0$. This proves that $d_r(f)=2$. So (b)
implies (a).\end{proof}

We now turn the characterisation of the equality $d_l(f)=2$. The following
corollary was first proved in \cite{C08} for irreducible morphisms in
standard components. Using \ref{prop:degree}, the proof given in
\cite{C08} generalizes to any Auslander-Reiten component. We thus
refer the reader to \cite{C08} for a detailed proof.
\begin{cor}
  \label{cor:degree_two}
Let $f\colon X\to Y$ be an irreducible
morphism with $X$ indecomposable. The following
conditions are equivalent:
\begin{enumerate}[(a)]
\item $d_l(f)=2$.
\item $Y$ the direct sum of at most two indecomposables;   $f$ is not minimal right
  almost split; and there exist $Z\in\ts{ind}\, A$ and a path of irreducible
  morphisms with composite $h$ lying on
  $\ts{rad}^2(Z,X)\backslash\ts{rad}^3(Z,X)$ and
  such that $fh=0$.
\item $\Gamma(\ts{mod}\, A)$ contains one of the two following
  configurations of meshes:
  \begin{align}
 &    
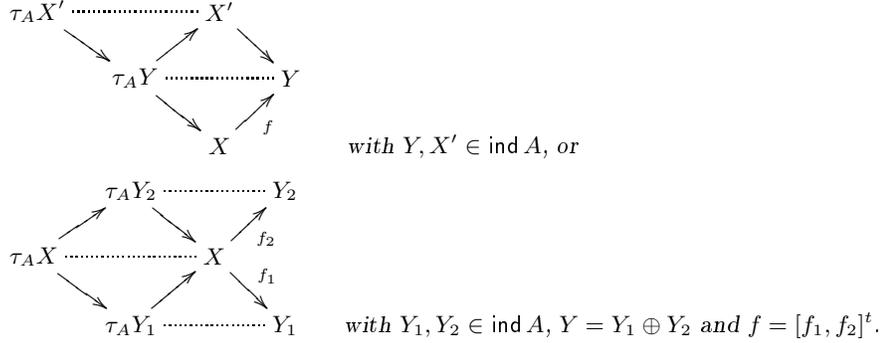
\xymatrix@=3ex{
\tau_A X' \ar@{.}[rr]\ar[rd] && X' \ar[rd]\\
&\tau_A Y \ar[ru] \ar[rd] \ar@{.}[rr] & & Y\\
&& X\ar[ru]_{f}& & \text{with $Y,X'\in\ts{ind}\,A$, or}
}\notag\\
&\xymatrix@=3ex{
& \tau_A Y_2 \ar[rd] \ar@{.}[rr] && Y_2\\
\tau_A X\ar[ru] \ar[rd] \ar@{.}[rr] & & X \ar[ru]_{f_2}
\ar[rd]^{f_1}\\
&\tau_A Y_1 \ar[ru] \ar@{.}[rr] && Y_1& \text{with
  $Y_1,Y_2\in\ts{ind}\,A$, $Y=Y_1\oplus Y_2$ and $f=[f_1,f_2]^t$.}
} \notag
  \end{align}
\end{enumerate}
\end{cor}

\section{Algebras of finite representation type}
\label{sec:s3}

In this section we prove our main theorem. First we need two lemmas.
\begin{lem}
\label{lem:simple}
  Let $S$ be a simple $A$-module, $S\hookrightarrow I$ its
  injective hull and $X\in\ts{ind}\, A$ such that $S$ is a direct
  summand of $\ts{soc}(X)$ . Assume that $I\twoheadrightarrow I/\ts{soc}(I)$ has
  finite left degree equal to $n$. Then
 there is a path in $\Gamma(\ts{mod}\,A)$ starting at $S$,
  ending at $I$, of length at most $n$ and going through $X$.
  In particular $X$, $S$  and $I$ lie in the same component of $\Gamma(\ts{mod}\, A)$.
\end{lem}
\begin{proof} 
Let $\Gamma$ be the component containing
$I$. We denote by $\pi$ the irreducible monomorphism $I\twoheadrightarrow
I/\ts{soc}(I)$ and by $\iota\colon S\hookrightarrow I$ the injective hull. It
follows from \ref{prop:kernel}  applied to  $\pi$ that $S\in\Gamma$ and
$\iota\in\ts{rad}^n(S,I)\backslash\ts{rad}^{n+1}(S,I)$.  Since $S$ is a direct
summand of $\ts{soc}(X)$, the injective hull $\iota$ factors through
$X$, that is, is equal to some composition
$S\xrightarrow{f}X\xrightarrow{g} I$.  Therefore 
there exist $l,m\geqslant 1$ such
that $f\in\ts{rad}^l(S,X)\backslash\ts{rad}^{l+1}(S,X)$,
$g\in\ts{rad}^m(X,I)\backslash\ts{rad}^{m+1}(X,I)$ and $l+m\leqslant
n$. Therefore $f$ (or $g$) is a sum of compositions of paths of
irreducible morphisms at least one of which has length $l$ (or $m$,
respectively). In particular, $P,X$ and $S$ all lie in $\Gamma$.
\end{proof}

Of course, \ref{lem:simple} has a dual statement which holds true
using dual arguments: If $S$ is a simple $A$-module with projective
cover $P\twoheadrightarrow S$ such that $\ts{rad}\,P\hookrightarrow P$ has
finite right degree equal to $n$ and if $S$ is a direct summand of $\ts{top}(X)$ for
some $X\in\ts{ind}\,A$, then there exists a path in
$\Gamma(\ts{mod}\,A)$ starting at $P$, ending at $S$, going through
$X$ and of length at most $n$. In particular, $P,X$ and $S$ all lie in
the same component of $\Gamma(\ts{mod}\,A)$.

\begin{lem}
\label{lem:unique_component}
  Assume $A$ is connected and that  for every indecomposable
  injective $I$ the quotient morphism
    $I\twoheadrightarrow I/\ts{soc}(I)$
has finite left degree. Let $n$ be the supremum of all these left
degrees. Then, for every $X\in\ts{ind}\,A$ there exists a path in
$\Gamma(\ts{mod}\,A)$ starting at $X$, ending at  some injective and
of length at most $n$. In particular, $\Gamma(\ts{mod}\,A)$ is finite
and connected.
\end{lem}
\begin{proof} The first assertion follows directly from
\ref{lem:simple}. In order to prove the second one it suffices to
prove that for all $I,J$  indecomposable injectives, $I$ and $J$ lie
on the same
component of $\Gamma(\ts{mod}\,A)$. Since $A$ is
connected there exists a sequence $I_0=I,I_1,\ldots,I_l=Q$ of
indecomposable injectives such that, letting $S_i=\ts{soc}(I_i)$, we
have that $S_i$ is a direct summand of
$\ts{soc}(I_{i-1}/\ts{soc}(I_{i-1}))$ or $S_{i-1}$ is a direct summand of
$\ts{soc}(I_i/\ts{soc}(I_i))$ for every $i=1,\ldots,l$. Accordingly,
\ref{lem:simple} implies that there exists a path $I_i\rightsquigarrow
I_{i-1}$ or $I_{i-1}\rightsquigarrow I_i$, respectively, in $\Gamma(\ts{mod}\,A)$.
This proves that $I$ and $J$ lie on the same component of
$\Gamma(\ts{mod}\,A)$. Since $\Gamma(\ts{mod}\,A)$ is locally finite,
this proves the lemma.
\end{proof}

Note that the dual statement of \ref{lem:unique_component} holds true
using dual arguments and the dual version of \ref{lem:simple}.
 In a previous version of this text, \ref{lem:unique_component} assumed
 an additional condition, dual to that on the degree of the morphisms
 $I\twoheadrightarrow I/\ts{soc}(I)$ for $I$ injective. The authors thank
 Juan Cappa for pointing out that this dual statement was unnecessary.

Now we can prove the main theorem. We recall its statement for
convenience.
\setcounter{Thm}{0}
\begin{Thm}
  Let $A$ be a connected finite dimensional $\ts k$-algebra over an
  algebraically closed field. The following conditions are equivalent:
  \begin{enumerate}[(a)]
  \item $A$ is of finite representation type. 
\item
For every indecomposable projective $A$-module $P$,
    the inclusion $\ts{rad}(P)\hookrightarrow P$ has finite right
    degree.
\item For every indecomposable injective $A$-module $I$, the quotient $I\to
  I/\ts{soc}(I)$ has finite left degree. 
\item For every irreducible epimorphism $f\colon X\to Y$ with $X$ or
  $Y$ indecomposable, the left degree of $f$ is finite.
\item For every irreducible monomorphism $f\colon X\to Y$ with $X$ or
  $Y$ indecomposable, the right degree of $f$ is finite.
  \end{enumerate}
\end{Thm}
\begin{proof}
 If $A$ is of finite representation type,
then $ \Gamma(\ts{mod}\, A)$ is connected  and $\ts{rad}^{\infty}=0$ (this follows
from the Lemma of Harada and Sai, for example) and the conditions (b) and (c) follow
from \ref{prop:kernel} and its dual. The implications $(b)\Rightarrow
(a)$ and $(c)\Rightarrow (a)$ follow from \ref{lem:unique_component}
and from its dual version, respectively. Thus, the conditions (a), (b)
and (c) are equivalent. Note that (d) implies (c), and (e) implies
(b). On the other hand, \ref{cor:degree_finite_type} and its dual
version show that (a) implies both (d) and (e). Therefore, the five
conditions (a), (b), (c), (d) and (e) are equivalent. 
\end{proof}

\begin{rem}
  Our arguments allow us to recover the following well-known
  implication using degrees of irreducible morphisms only:
\emph{ If $\ts{rad}^{\infty}=0$ then $A$ is of finite representation type
    and $\Gamma(\ts{mod}\,A)$ is connected}.
Indeed, if $\ts{rad}^{\infty}=0$ then both (b) and (c) hold true in
Theorem~\ref{thm1}. So $A$ is of finite representation type and
$\Gamma(\ts{mod}\,A)$ is connected.
\end{rem}

\section{Composition of morphisms}
\label{sec:s4}

Let $A$ be a finite dimensional $\ts k$-algebra and $\Gamma$ be a  component
of $\Gamma(\ts{mod}\, A)$. In view of \ref{rem:domain_explode}, (b), there seems
to be a connection between the degree of an irreducible morphism and
the behavior of the composite of $n$ irreducible morphisms between
indecomposable modules (for any $n$). This motivates the work of the
present section, that is, to study when the composite of $n$
irreducible morphisms between indecomposable modules lies in
$\ts{rad}^{n+1}$. The following result characterizes such a situation
when $\Gamma$ has \emph{trivial valuation} (that is, has no multiple arrows).
\begin{prop}
\label{prop:composite}
Let $n\geqslant 1$ be an
  integer and $X_1,\ldots,X_{n+1}$ be modules in $\Gamma$. Consider the
  following assertions:
  \begin{enumerate}[(a)]
  \item There exist irreducible morphisms $h_i\colon X_i\to
    X_{i+1}$ for every $i$ such that $h_n\ldots
    h_1\in\ts{rad}^{n+1}(X_1,X_{n+1})\backslash\{0\}$.
  \item There exist irreducible morphisms $f_i\colon X_i\to
    X_{i+1}$ together with morphisms $\varepsilon_i\colon X_i\to
    X_{i+1}$ such that $f_n\ldots f_1=0$,
    $\varepsilon_n\ldots \varepsilon_1\neq 0$ and $\varepsilon_i=f_i$ or
    $\varepsilon_i\in\ts{rad}^2(X_i,X_{i+1})$ for every $i$.
  \end{enumerate}
Then (b) implies (a). Also, if $h_1,\ldots,h_n$ satisfy (a) and
represent arrows with trivial valuation, then (b) holds true. In
particular, (a) and (b) are equivalent if $\Gamma$ has trivial valuation.
\end{prop}
\begin{proof} 
 Let
$h_i\colon X_i\to X_{i+1}$, $i=1,\ldots,n$, be irreducible
morphisms in $\ts{ind}\, A$
such that $h_n\ldots
h_1\in\ts{rad}^{n+1}(X_1,X_{n+1})\backslash\{0\}$ and such that the
arrows represented by $h_1,\ldots,h_n$ have trivial valuation. Let $F\colon
\ts k(\widetilde{\Gamma})\to\ts{ind}\,\widetilde{\Gamma}$ be a well-behaved
functor with respect to the generic covering
$\pi\colon\widetilde{\Gamma}\to\Gamma$ and let $x_1\in F^{-1}(X_1)$. Since $\pi\colon \widetilde{\Gamma}\to \Gamma$ is a
covering of quivers and the arrow represented by $h_1$ has trivial
valuation, there is exactly
one arrow $x_1\xrightarrow{\alpha_1}x_2$ in
$\widetilde{\Gamma}$ starting from $x_1$ and such that $Fx_2=X_2$. By repeating the
same argument, we deduce that there is exactly one path in $\widetilde{\Gamma}$:
\begin{equation}
  x_1\xrightarrow{\alpha_1}x_2\to\ldots\to x_n\xrightarrow{\alpha_n}x_{n+1}\notag
\end{equation}
starting from $x_1$ of length $n$ and such that $Fx_i=X_i$ for every
$i$. Let $i\in\{1,\ldots,n\}$, then $F(\overline\alpha_i)\colon X_i\to
X_{i+1}$ is irreducible so that $h_i=\lambda_iF(\overline\alpha_i)+h_i'$,
where $\lambda_i\in \ts k^*$ and $h_i'\in\ts{rad}^2(X_i,X_{i+1})$ because
$\pi(\alpha_i)$, represented by $h_i$, has trivial valuation.  Since $h_n\ldots
h_1\neq 0$, we have a non-zero morphism:
\begin{equation}
\begin{array}{l}
  \lambda F(\overline{\alpha_n\ldots\alpha_1})\\ +
  \sum\limits_{t=1}^n\sum\limits_{i_1<\ldots<i_t} F(\overline\alpha_n) \ldots
  F(\overline\alpha_{i_t+1}) h_{i_t}' F(\overline\alpha_{i_t-1})\ldots
  F(\overline\alpha_{i_1+1}) h_{i_1}' F(\overline\alpha_{i_1-1})\ldots F(\overline\alpha_1)
\end{array}\tag{$\star$}
\end{equation}
where $\lambda=\lambda_1\ldots\lambda_n\in \ts k^*$ and the whole sum lies
on $\ts{rad}^{n+1}(X_1,X_{n+1})$. In particular,
$F(\overline{\alpha_n\ldots\alpha_1})$ lies on $\ts{rad}^{n+1}(X_1,X_{n+1})$. By
Theorem~\ref{prop:covering}, we have
$\overline{\alpha_n\ldots\alpha_1}\in\mathfrak{R}^{n+1}\ts k(\widetilde{\Gamma})(x_1,x_{n+1})$. Since $\widetilde{\Gamma}$
is a component with length, we deduce that
$\overline{\alpha_n\ldots\alpha_1}=0$ and therefore $F(\overline\alpha_n)\ldots
F(\overline\alpha_1)=0$. This and ($\star$) imply that there exist
$t\in\{1,\ldots,n\}$ and $i_1<\ldots<i_t$ such that:
\begin{equation}
  F(\overline\alpha_n) \ldots
  F(\overline\alpha_{i_t+1}) h_{i_t}' F(\overline\alpha_{i_t-1})\ldots
  F(\overline\alpha_{i_1+1}) h_{i_1}' F(\overline\alpha_{i_1-1})\ldots
  F(\overline\alpha_1)\neq 0\ .\tag{$\star\star$}
\end{equation}
We thus let:
\begin{enumerate}
\item $f_i=F(\overline\alpha_i)$ for every $i\in\{1,\ldots,n\}$. So
  $f_i\colon X_i\to X_{i+1}$ is irreducible because $F\colon \ts k(\widetilde{\Gamma})\to
  \ts{ind}\,\Gamma$ is well-behaved.
\item $\varepsilon_{i_j}=h_{i_j}'$ for every $j\in\{1,\ldots,t\}$. So
  $\varepsilon_{i_j}\in \ts{rad}^2(X_{i_j},X_{i_j+1})$.
\item $\varepsilon_i=f_i$ for every $i\in\{1,\ldots,n\}\backslash\{i_1,\ldots,i_t\}$.
\end{enumerate}
In particular, $\varepsilon_n\ldots\varepsilon_1\neq 0$ because of ($\star\star$).
The morphisms $f_i$ and $\varepsilon_i$ ($i\in\{1,\ldots,n\}$) satisfy the
conclusion of (b). This proves (b) when $h_1,\ldots,h_n$ satisfy (a)
and represent arrows with trivial valuation.

For the implication (b) implies (a), we refer the reader to the proof
of \cite[Thm. 2.7]{CT08} (where the standard hypothesis made therein
is not used for that implication).

The equivalence between (a) and (b) when $\Gamma$ has trivial valuation
follows from the above considerations.
\end{proof}

\begin{rem}
  Let $h_1,\ldots,h_n$, $i=1,\ldots, n$, be morphisms satisfying (a) in
  \ref{prop:composite}. Under additional assumption such as,
  $\alpha(\Gamma)\leqslant 2$ (\cite{CPT04}), or $n=2$ (\cite{CCT07}), or
  $n=3$ (\cite{CCT08a}) or the path $h_1,\ldots,h_n$ is almost
  sectional (\cite{CCT08}), it is known that the arrows in $\Gamma$
  represented by $h_1,\ldots,h_n$ all have trivial valuation. However,
  it is still an open question to know whether this is always the case.
\end{rem}

Our last result concerns sums of composites of paths in a sectional
family (\ref{defn:sectional_family}). Note that this result
extends the well-known result of Igusa and Todorov (\cite{IT84})
and which asserts that if
$\cdot\xrightarrow{f_1}\cdot\to\cdots\to\cdot\xrightarrow{f_l}\cdot$
is a sectional path of irreducible morphisms between indecomposables,
then the composite $f_n\cdots f_1$ does not lie in $\ts{rad}^{n+1}$
and, therefore, is non-zero. Recall that a sectional path of
irreducible morphisms between indecomposables is a particular case of
a sectional family of paths (\ref{rem:sectional_family}).
\begin{prop}
Let $X,Y$ be indecomposable modules in $\Gamma$. Let
$\{X\xrightarrow{f_{i,1}}X_{i,1}\to\cdots\to
  X_{i,l_{i-1}}\xrightarrow{f_{i,l_i}} Y\}_{i=1,\ldots,r}$ be
  a sectional family of paths starting in $X$, ending in $Y$ and of irreducible morphisms between
  indecomposables. Let
  $n=\underset{i=1,\ldots,r}{\ts{min}}l_i$. Then $\sum\limits_{i=1}^rf_{i,l_i}\cdots
  f_{i,1}$ lies in $\ts{rad}^n(X,Y)$ and does not lie in
  $\ts{rad}^{n+1}(X,Y)$. In particular it is non-zero.
\end{prop}
\begin{proof}
Let $\pi\colon \widetilde{\Gamma}\to\Gamma$ be the generic covering
and let $x\in\pi^{-1}(X)$. We apply \ref{prop:well_behaved_sectional}
from which we adopt the
notations ($x_{i,j},\alpha_{i,j}$). In particular, there exists a
well-behaved functor $F\colon
\ts k(\widetilde{\Gamma})\to \ts{ind}\,\Gamma$ such that
$F(\overline{\alpha_{i,j}})=f_{i,j}$ for every $i,j$.

For every $y\in \pi^{-1}(Y)$ let $I_y$ be the set of indices such that
$x_{i,l_i}=y$. For each $i$, let $u_i$ be the path $x\xrightarrow{\alpha_{i,1}}x_{i,1}\to\cdots\to
x_{i,l_i-1}\xrightarrow{\alpha_{i,l_i}} x_{i,l_i}$. Therefore, there
exists some $y_0\in \pi^{-1}(Y)$ such
that $I_{y_0}$ is non-empty and all the paths  $u_i$, $i\in I_{y_0}$,
have length $n$. Moreover, each path $u_i$, for $i\in\{1,\ldots,r\}$
is sectional, because $\pi\colon\widetilde{\Gamma}\to\Gamma$ is a
covering of translation quivers and $\pi(u_i)$ is a sectional path
$X\to X_{i,1}\to \cdots\to X_{i,l_i-1}\to X_{i,l_i}$
(\ref{rem:sectional_family}, (1)).

Clearly, the sum $\sum\limits_{i=1}^rf_{i,l_i}\cdots
  f_{i,1}$ equals $\sum\limits_{i=1}^rF(\overline{u_i})$ and lies in
  $\ts{rad}^n(X,Y)$. By absurd, assume that it lies
  in $\ts{rad}^{n+1}(X,Y)$. Using Theorem~\ref{prop:covering}, we deduce that
  $\sum\limits_{i\in
    I_y}\overline{u_i}\in\mathfrak{R}^{n+1}\ts k(\widetilde{\Gamma})(x,y)$,
  for every $y\in\pi^{-1}(Y)$. In particular, $\sum\limits_{i\in
    I_{y_0}}\overline{u_i}\in\mathfrak{R}^{n+1}\ts k(\widetilde{\Gamma})(x,y_0)$. This
  contradicts \ref{prop:generic_covering}, (e), because the paths
  $u_i$, for $i\in I_{y_0}$, are sectional and of length $n$.
\end{proof}

\begin{acknowledgements}
The authors would like to thank Shiping Liu for fruitful discussions
 which led to improvements in the text. 
This text was elaborated during a stay of the second author at Mar del
Plata, he would like to thank Claudia Chaio, Sonia Trepode and all the
algebra research group of Mar del Plata for their warm hospitality.
Finally all authors would like to
 thank Ibrahim Assem for his encouragements and constant interest in
 this work.
\end{acknowledgements}

\bibliographystyle{plain}
\bibliography{biblio}

\affiliationone{
Claudia Chaio and Sonia Trepode\\
Departamento de Matem\'atica, FCEyN, Universidad
Nacional de Mar del Plata, Funes 3350, 7600 Mar del Plata, Argentina\\
\email{algonzal@mdp.edu.ar\\
strepode@mdp.edu.ar}}
\affiliationtwo{
Patrick Le Meur\\
CMLA, ENS Cachan, CNRS, UniverSud, 61 Avenue du President Wilson,
F-94230 Cachan\\
\email{
Patrick.LeMeur@cmla.ens-cachan.fr}}
\end{document}